\documentclass[12pt]{amsart}
\newtheorem{thm}{Theorem}[subsection]
\newtheorem{cor}[thm]{Corollary}
\newtheorem{lem}[thm]{Lemma}

\theoremstyle{remark}

\theoremstyle{remark}

\numberwithin{equation}{subsection}

\usepackage{bbold}
\usepackage[T2A]{fontenc}

\usepackage{CJKutf8}

\usepackage{amsmath,amstext,amsfonts,amsthm,amssymb}
\usepackage{eucal}
\usepackage{graphicx}
\usepackage{tikz-cd}
\usetikzlibrary{patterns}
\usetikzlibrary{calc}

\swapnumbers

{  \theoremstyle{definition}

}
\newcommand{\lra}{\longrightarrow}

\newcommand{\bea}{\begin{eqnarray*}}
\newcommand{\eea}{\end{eqnarray*}}
\newcommand{\bean}{\begin{eqnarray}}
\newcommand{\eean}{\end{eqnarray}}

\newcommand{\xni}{x^{(n)}_i}
\newcommand{\partxni}{\partial^{(n)}_{x_i}}
\newcommand{\dxni}{dx^{(n)}_i}
\newcommand{\partdxni}{\partial_{dx_i}^{(n)}}

\newcommand{\fg}{\mathfrak g}

\newcommand{\bun}{\mathbb1}

\newcommand{\CF}{\mathcal{F}}

\newcommand{\CH}{\mathcal{H}}

\newcommand{\CO}{\mathcal{O}}

\newcommand{\CT}{\mathcal{T}}

\newcommand{\BC}{\mathbb{C}}

\newcommand{\BZ}{\mathbb{Z}}

\newcommand{\nc}{\newcommand}

\nc{\Id}{\text{Id}}
\nc{\la}{\lambda}

\begin{document}

\title{ chiral cartan calculus}

\author{ Fyodor Malikov, Vadim Schechtman, Boris Tsygan}

\thanks{The results of this paper have been reported by the second named author on November, 2023  at SISSA, Trieste. He is grateful to the organizers for the invitation.}

\thanks{}

\maketitle

\begin{center}
 {\it To the blessed memory of Yuri Ivanovich Manin}
\end{center}













 
   \quad\quad

\section{introduction}
\label{introduction}
\

\

The aim of this note is a study of the Gauss - Manin connection on the chiral de Rham complex.

In the seminal paper \cite{G} Grothendieck has explained (and introduced the name of) the Gauss - Manin connection, cf. \cite{M}; his explanation is based on the Cartan homotopy 
formula ("formule chère à Cartan"). In \cite{T} (see also \cite{GDT}, \cite{Ge}, \cite{DTT}) the mechanism of the Cartan homotopy formula was formalized under the name 
of {\it calculus} and generalized to the noncommutative situation. 

A calculus is a structure that formalizes the standard relations between the graded Lie algebra of polyvector fields 
$\Lambda^\bullet\CT_X$ and the de Rham complex $\Omega^\bullet_X$ of a manifold $X$.  

The aim of the present note is to propose a chiral version of the calculus (we deal here only with the commutative case, although we believe that a noncommutative version exists as well).

In sect.~\ref{calculus} we propose our version of the classical formalism of Grothendieck and 
Katz - Oda \cite{KO}.
It is a more precise version of the Grothendieck construction. 

In  sect.~\ref{the relative chiral de rham complex} we discuss the relative chiral de Rham complex $\Omega_{X/S}^{ch}$ for a smooth morphism 
$f: X\lra S$ of smooth algebraic varieties, the usual chiral de Rham complex [MSV] corresponding to the case $S = *$. In sect.~\ref{chiral calculus} we introduce a chiral version of the calculus algebra which acts on 
$\Omega_{X/S}^{ch}$. Here a mechanism of anomaly cancellation appears, similar to 
the one which gives rise to the usual chiral de Rham. 

After passing to the cohomology this action 
induces the usual Gauss - Manin connection on  the relative cohomology sheaf $R^\bullet f_*\CO_X$.

\bigskip

The main results of this paper are stated and proved in sects.~\ref{calculus} and  \ref{chiral calculus}, and the reader might find it reasonable to skip sect.~\ref{the relative chiral de rham complex} on first reading. Sect.~\ref{the relative chiral de rham complex} is our attempt, perhaps not overly successful, to retell a number of well-known results using the language of chiral algebras introduced in \cite{BD}. The reason we felt this was worth doing is best understood by contemplating the exact sequence (\ref{chiral picard lie}). This exact sequence represents the Beilinson-Drinfeld way of defining an algebra of chiral differential operators and the degree of clarity it  achieves is quite extraordinary, but it does  use  the technology of \cite{BD} in earnest -- compound pseudo-tensor categories, Lie* algebras and algebroids, etc.,  the concepts which do not lend themselves easily to treatment by more traditional vertex algebra methods.

\section{calculus}
\label{calculus}

\subsection{Calc 1.}
\label{Calc1.}
 Let $p: X\rightarrow S$ be a smooth morphism of smooth algebraic varieties over $\BC$.  We will be writing $X/S$ to indicate that this morphism is fixed
 and various standard constructions are understood in the relative sense. For example, the structure sheaf $\CO_{X/S}$ is automatically a sheaf of $p^{-1}\CO_S$-modules, $\Omega_{X/S}$ will stand for the module of K\"ahler differentials over $\CO_S$, $\CT_{X/S}$ is the algebra of $\CO_S$-derivations, etc.
 
 Along with the relative tangent bundle $TX/S$, one can consider  its superversion, $\Pi TX/S$, obtained from the former  declaring that the fibers of 
 the projection $TX/S\rightarrow X/S$ to be odd.
 
The supervariety $\Pi TX/S$ carries an algebra of differential operators, $D_{\Pi TX/S}$, which is filtered, $D_{\Pi TX/S}=\cup_{i\geq 0} D_{\Pi TX/S}^{\leq i}$, and generated by the standard Picard-Lie algebroid:
\begin{equation}
\label{ordinary picard lie}
0\rightarrow \Omega^\bullet(X/S)\rightarrow D_{\Pi TX/S}^{\leq 1}\rightarrow \CT_{\Pi TX/S}\rightarrow 0,
\end{equation}
where the de Rham complex $\Omega^\bullet(X/S)$, regarded for the time being only as a supercommutative algebra, is  nothing but the structure
 sheaf $\CO_{\Pi TX/S}$. The canonical action of vector fields on differential forms defines  a Lie algebra embedding
\begin{equation}
\label{formula for lie br on diffe forms}
\mbox{Lie}: \CT_{X/S}\rightarrow D_{\Pi TX/S}^{\leq 1},\; \xi\mapsto \mbox{Lie}_{\xi}.
\end{equation}

Similarly, we have a  morphism of Lie algebra modules
\begin{equation}
\label{formula for contr with a vect fie on diffe forms}
\iota: \Pi\CT_{X/S}\rightarrow D_{\Pi TX/S}^{\leq 1},\;,
\end{equation}
so that
\begin{equation}
\label{class form 1}
[\mbox{Lie}_\xi,\iota_\eta]=\iota_{[\xi,\eta]},\; [\iota_{\xi},\iota_\eta]=0\; \forall \xi,\eta\in\CT_{X/S}.
\end{equation}

We can choose, locally on $X/S$, what in the D-module theory  is habitually called a coordinate system, that is, a collection $\{x_i,\partial_{x_i}, 1\leq i\leq \mbox{dim}\CT_{X/S}\}$,  where $\{\partial_{x_i}\}$ is a local basis of $\CT_{X/S}$, $\{x_i\}\subset \CO_X$ and $\partial_{x_i}(x_j)=\delta_{ij}$.  We can further extend it to  a coordinate system on $\Pi TX/S$:  $\{x_i, dx_i, \partial_{x_i}, \partial_{dx_i}\}$, where the differential $dx_i$ is regarded as an odd function on $\Pi TX/S$ linear along the fibers of the projection $\Pi TX/S\rightarrow X/S$, and $ \partial_{dx_i}$ is an odd derivation such that 
$ \partial_{dx_i}(dx_j)=\delta_{ij}$ and $ \partial_{dx_i}(x_j)=0$. Given this choice, we
can write explicit formulas, which will prove important in the chiral algebra situation:
\begin{equation}
\label{expl formula for lie derib and contra ordin}
\mbox{Lie}_{f_i\partial_{x_i}}=f_i\partial_{x_i}+\frac{\partial f_i}{\partial x_\alpha}dx_\alpha\partial_{dx_i},\;  \iota_{f_i\partial_{x_i}}=f_i\partial_{dx_i}.
\end{equation}
{\em Here and elsewhere, the summation w.r.t. repeated indices will be tacitly assumed.}

There is another celebrated element of  $\CT_{\Pi TX/S}$, the de Rham differential, $d_{X/S}=dx_i\partial_{x_i}$, and it satisfies
\begin{equation}
\label{comm rel invo de rham diff}
[d_{X/S},d_{X/S}]=0,\; [d_{X/S},\iota_\xi]=\mbox{Lie}_\xi,\; [d_{X/S},\mbox{Lie}_\xi]=0.
\end{equation}

To conclude:

\begin{lem}
\label{summ of commu rel on the ordn deram}
The maps (\ref{formula for lie br on diffe forms}, \ref{formula for contr with a vect fie on diffe forms}) define a morphism of differential graded Lie algebras
\begin{equation*}
\Pi\CT_{X/S}\rtimes\CT_{X/S}\rightarrow D_{\Pi TX/S}^{\leq 1},
\end{equation*}
the differential on the right hand side being $[d_{X/S},.]$ .
\end{lem}

\subsection{ }
\label{boryas constr ordin}
Given the morphism $p: X\rightarrow S$, sect. \ref{Calc1.}, consider a connection 
\begin{equation}
\label{def of cconn ordin}
\nabla:
\CO_{X/S}\rightarrow \Omega_S\otimes_{\CO_S}\CO_{X/S}
\end{equation}
such that
\begin{equation}
\label{def prop conn ordin}
\nabla(fg)=f\nabla(g)+g\nabla(f), \nabla(h)=d_Sh\;\forall f.g\in\CO_{X/S}, h\in\CO_S,
\end{equation}
where $d_Sh$ is the differential of $h$. Requirements (\ref{def prop conn ordin}) say that, geometrically, $\nabla$ is a rule that lifts a vector field on $S$ to that on $X/S$. Such connections exist at least locally on $X$.

Given $\nabla$, we obtain its curvature
\begin{equation}
\label{def of curv ordin}
R=\frac{1}{2}[\nabla,\nabla].
\end{equation}
Choosing, locally, a coordinate system, $\{u_i,\partial_{u_i}\}$, on $S$, we obtain:
\begin{equation*}
\nabla=du_i(\partial_{u_i}+\xi_i),\; R=\frac{1}{2}du_idu_j(\partial_{u_i}\xi_j-\partial_{u_j}\xi_i+[\xi_i,\xi_j]),\;\xi_i,\xi_j\in\CT_{X/S}.
\end{equation*}

We will now extend the scalars from $\CO_S$ to the entire de Rham complex on $S$, $\Omega^\bullet_S$, thus replacing $\Omega^\bullet_{X/S}$ with
$\Omega^\bullet_S\otimes_{\CO_S}\Omega^\bullet_{X/S}$.  Just as given $\xi\in\CT_{X/S}$, $\mbox{Lie}_\xi$ defines a derivation of 
$\Omega^\bullet_{X/S}$,  the connection and curvature can be promoted to derivations of $\Omega^\bullet_S\otimes_{\CO_S}\Omega^\bullet_{X/S}$, to be denoted $\mbox{Lie}_\nabla$ and $\mbox{Lie}_R$ resp.; e.g. in local coordinates $\mbox{Lie}_\nabla=du_i(\partial_{u_i}+\mbox{Lie}_{\xi_i})$. The morphism $\iota$, see (\ref{formula for contr with a vect fie on diffe forms}), affords a similar extension, and we obtain derivations such as $\iota_R$. 

Consider now the operator, in fact a derivation, 
\begin{equation*}
D_\nabla=\mbox{Lie}_\nabla+d_{X/S}-\iota_R\in\mbox{End}_\BC(\Omega^\bullet_S\otimes_{\CO_S}\Omega^\bullet_{X/S})
\end{equation*}

\begin{lem}
\label{boryas diff ordina}
The pair $(\Omega^\bullet_S\otimes_{\CO_S}\Omega^\bullet_{X/S}, D_\nabla)$ is a complex, i.e.,
\[
D_\nabla^2=0.
\]
\end{lem}

 {\em Proof:}
 
 The result follows from the following computation:
  \begin{eqnarray*}
  & &2 D_\nabla^2=[ D_\nabla,  D_\nabla]=\\
 & & [\mbox{Lie}_\nabla,\mbox{Lie}_\nabla]+[ d_{X/S}, d_{X/S}]+[\iota^{ch}_{R},\iota^{ch}_{R}]\\
 &+&2[\mbox{Lie}_\nabla,d_{X/S}]-2[\mbox{Lie}_\nabla, \iota_{R}]-2[d_{X/S},\iota_{R}]=\\
 & &2 \mbox{Lie}_{R}-2\iota_{[\nabla,R]}-2\mbox{Lie}_{R}=0,
 \end{eqnarray*}
 which consists in a repeated application of the `calculus' of sect. \ref{Calc1.}, that is, classical differential-geometric formulas (\ref{class form 1}, \ref{comm rel invo de rham diff}); one exception is the term $\iota_{[\nabla,R]}$, which occurs thanks to $\iota_{[\nabla,R]}=[\mbox{Lie}_\nabla, \iota_{R}]$ (formula
(\ref{class form 1})), and it vanishes thanks to the Bianchi identity $[\nabla,R]=0$. $\qed$

We will call the complex $(\Omega^\bullet_S\otimes_{\CO_S}\Omega^\bullet_{X/S}, D_\nabla)$ the  
{\it local connection complex}, and denote it $CC(X/S, \nabla)$.  It is an example of a {\it twisted complex} in the sense of Toledo - Tong, cf. \cite{TT}.  

The differential $D_\nabla$ is not unique, as it depends on a choice of $\nabla$, but various choices give rise to canonically isomorphic complexes.

 \begin{lem}
  \label{iso of diff connectns ordinary}
  The map
  \[
  \exp{(\iota_{\nabla_2-\nabla_1})}:(\Omega^\bullet_S\otimes_{\CO_S}\Omega^{\bullet}_{X/S},D_{\nabla_1})\rightarrow
  (\Omega^\bullet_S\otimes_{\CO_S}\Omega^{\bullet}_{X/S},D_{\nabla_2})
  \]
  is an isomorphism of complexes.
  \end{lem}
  
  {\em Sketch of Proof:}
  
  It is clear that the operator $\iota_{\nabla_2-\nabla_1}\in\mbox{End}_\BC(\Omega^\bullet_S\otimes_{\CO_S}\Omega^{\bullet}_{X/S})$ is nilpotent
  and is an even derivation of the algebra $\Omega^\bullet_S\otimes_{\CO_S}\Omega^\bullet_{X/S}$. Therefore it defines an automorphism of 
  $\Omega^\bullet_S\otimes_{\CO_S}\Omega^\bullet_{X/S}$.   What we need to check then is that 
  \[
   \exp{(ad_{\iota_{\nabla_2-\nabla_1,(0)}})}(D_{\nabla_1})=D_{\nabla_2},
   \]
   where $ad_{\iota_{\nabla_2-\nabla_1}}$ means the operation of  taking  the bracket, $[\iota_{\nabla_2-\nabla_1},.]$. A direct  computation, using the various identities as in the proof of  Lemma \ref{iso of diff connectns ordinary}, will show that
   \begin{eqnarray*}
   \exp{(ad_{\iota_{\nabla_2-\nabla_1}})}(D_{\nabla_1})&=&D_{\nabla_1}+[\iota_{\nabla_2-\nabla_1}, D_{\nabla_1}]\\
   &+&
   \frac{1}{2}[\iota_{\nabla_2-\nabla_1},[\iota_{\nabla_2-\nabla_1}, D_{\nabla_1}]]=D_{\nabla_2}.\; \qed
   \end{eqnarray*}
   
   \subsection{Comparison with Katz-Oda.}
   \label{Comparison with Katz-Oda.}
   The complex $(\Omega^\bullet_S\otimes_{\CO_S}\Omega^\bullet_{X/S}, D_\nabla)$ is filtered, 
   $$
   F^p(\Omega^\bullet_S\otimes_{\CO_S}\Omega^\bullet_{X/S})=\Omega^{\geq p}_S\otimes_{\CO_S}\Omega^\bullet_{X/S}
   $$ 
   and the first term of the corresponding spectral sequence
   $$
   E_1:\ \Omega^\bullet_S\otimes H^\bullet_{DR}(X/S,\nabla) \overset{\nabla_{GM}}\lra 
   \Omega^{\bullet + 1}_S\otimes H^\bullet_{DR}(X/S,\nabla)
   $$

obviously coincides with that of $\Omega^\bullet_X$ introduced in \cite{KO}, for any $\nabla$, thus reproducing the Katz-Oda construction \cite{KO} of the Gauss-Manin connection. Cf. also \cite{SV}.
   
   \subsection{Globalization.}
   \label{Globalization-Ordin.}
   We will now show that Lemmas~\ref{boryas diff ordina} and~\ref{iso of diff connectns ordinary} imply the existence of a globally defined complex; such a complex is not unique,   but different such  complexes  are canonically isomorphic. 
   
   The variety $S$ being smooth, it can be covered by open subsets $S_\alpha$ \'etale over $\BC^n$, hence each having a coordinate system $\{u_i,\partial_{u_i}\}$, $\partial_{u_i}(u_j)=\delta_{ij}$. The morphism $p:X\rightarrow S$ being smooth, the preimage $p^{-1}(S_\alpha)$ has an open cover $\{U_{\alpha\beta}\}$ so that each $U_{\alpha\beta}$ has a coordinate system $\{\{u_i,x_j, \partial_{u_i}, \partial_{x_j}\}$. This means that over each  $U_{\alpha\beta}$, the restriction
 $\Omega^\bullet_S\otimes_{\CO_S}\Omega^{\bullet}_{X/S}|_{U_{\alpha\beta}}$ regarded as a sheaf over $S_\alpha$ carries a connection, 
 $\nabla_{\alpha\beta}$.   
 
 To streamline the notation: we have a sheaf $\Omega\stackrel{\mbox{def}}{=}\Omega^\bullet_S\otimes_{\CO_S}\Omega^{\bullet}_{X/S}$ on $X$ and an open
 cover $\{U_i\}$ of $X$;  on each restriction $\Omega_i\stackrel{\mbox{def}}{=}\Omega|_{U_i}$ we have chosen a connection $\nabla_i$ hence a differential $D_{\nabla_i}:\Omega_{i}\rightarrow \Omega_{i}$, Lemma~\ref{boryas diff ordina}.  According to Lemma~\ref{iso of diff connectns ordinary}, over the intersection $U_i\cap U_j$ there is an
 isomorphism 
 \begin{equation*}
  G_{ji}\stackrel{\mbox{def}}{=}\exp{(\iota_{\nabla_j-\nabla_i})}: (\Omega_i,D_{\nabla_i})\rightarrow (\Omega_j,D_{\nabla_j}).
  \end{equation*}
  It is clear that on triple intersections we have the consistency condition:
  \begin{equation*}
  G_{ki}= G_{kj}\circ G_{ji}.
  \end{equation*}
  Therefore, having fixed a cover $\{U_i\}$ and a bunch of connections $\{\nabla_i\}$, which we will denote simply by $\nabla$, we obtain a new globally defined complex, $\Omega_{\nabla}$,  by tearing the original 
  $\Omega^\bullet_S\otimes_{\CO_S}\Omega^{\bullet}_{X/S}$ apart and regluing the pieces $\{(\Omega_i,\nabla_i)\}$ using  $\{G_{ij}\}$ as transition functions. This proves the existence.
  
  Two distinct such complexes, $\Omega_{\nabla}$ and $\Omega_{\nabla'}$ (attached to distinct collections $\nabla$ and $\nabla'$) are canonically isomorphic. To see this, assume, as we can, that both are defined using the same open cover $\{U_i\}$. Then the collection of isomorphisms 
  \begin{equation*}
  \exp{(\iota_{\nabla'_i}-\iota_{\nabla_i})}:\Omega_{\nabla}|_{U_i}\rightarrow\Omega_{\nabla'}|_{U_i}
  \end{equation*}
  well defines an isomorphism of complexes $\Omega_{\nabla}\rightarrow\Omega_{\nabla'}$. 
  
  We will call the resulting glued complex defined up to a canonical isomorphism the {\it global connection complex}, and denote it $CC(X/S)$; it gives rise to the Gauss - Manin connection on the relative cohomology.

\section{the relative chiral de rham complex}
\label{the relative chiral de rham complex}

\subsection{ }
\label {d-mod-gener}

Let $C$ be a smooth algebraic curve over $\BC$, $D_C$ the algebra of differential operators on $C$, $\mbox{Mod}^l(D_C)$, $\mbox{Mod}^r(D_C)$  the category of left (right, resp.) $D_C$-modules.  

For a finite set $I$, we let $C^I$ be the Cartesian product, and for a surjection $\pi^{(J/I)}:J\rightarrow I$ we denote by $\Delta^{(J/I)}$ the diagonal map $C^I\hookrightarrow C^J$; if $I$ is a point, then we write $\Delta^{(J)}: C\hookrightarrow C^J$.

Given an $I$-family of $D_C$-modules, $\{M_i\}$, left or right, we have a $D_{C^I}$-module, left or right resp., $\boxtimes_IM_i$, and if $\{M_i\}\subset \mbox{Mod}^r(D_C)$ we obtain
 \[
\Delta^{(J/I)}_\ast(\boxtimes_IM_i)\in \mbox{Mod}^r(D_{C^J}),
\]
 where $\Delta^{(J/I)}_*$ is the standard direct image functor: 
 \[
 \Delta^{(J/I)}_\ast(\boxtimes_IM_i)=(\boxtimes_IM_i)\otimes_{D_{C^I}}D_{C^I\rightarrow C^J}.
 \]
 The $D_{C^J}$-module structure of $ \Delta^{(J/I)}_\ast(\boxtimes_IM_i)$ becomes apparent  locally  where we can include
  a local coordinate system on $C^I$ into that on $C^J$. For example, if $C=\BC$ with coordinate $z$ and, for simplicity, $I$ is a point, then for any $j_0\in J$
 \begin{equation}
 \nonumber
 \Delta^{(J)}_\ast(M)\stackrel{\sim}{\rightarrow}M[\partial_{z_j},\; j\neq j_0],
 \end{equation}
 and the action of $D_{C^J}$ is determined by the following requirements: $\CO_{C^{J}}$ acts on $M$ via the projection 
 $\BC[z_j, j\in J]\rightarrow \BC[z]$, $z_j\mapsto z$, the action of $\partial_j$, $j\neq j_0$ is obvious, and $\partial_{z_{j_0}}$ acts as
 $\partial_z-\sum_{j\neq j_0}\partial_{z_j}$; of course, the action of $\partial_z$ and $\BC[z]$ on $M$ is the tautological one.
 
 This means that $ \Delta^{(J)}_\ast(M)$ is the space of maps  $\CO_{C^J}\rightarrow  \Delta^{(J)}_\bullet$ that are differential operators in the directions transversal to the diagonal which suggests a convenient "functional" realization so popular in the vertex algebra theory, e.g.
 \begin{equation}
 \label{functional-realiz}
 \Delta^{(J)}_\ast(M)\stackrel{\sim}{\rightarrow}M[(z_j-z_{j_0})^{-1},\; j\neq j_0]\bigwedge_{j\neq j_0}\frac{dz_j}{z_j-z_{j_0}}.
 \end{equation}
 The element of $ \Delta^{(J)}_\ast(M)$ that is represented by $m\otimes \prod_{j\neq j_0}\partial_{z_j}^{k_j}/k_j!$ or, equivalently,  by
 $m\otimes \bigwedge_{j\neq j_0}dz_j/(z_j-z_{j_0})^{-k_j-1}$ maps
 \begin{equation}
 \nonumber
 f\mapsto m\cdot \prod_{j\neq j_0}\partial_{z_j}^{k_j}/k_j!f|_{z_j=z_{j_0}}
 \end{equation}
 Of course, there is nothing special in the choice of coordinates $\{z_j-z_{j_0}, z_{j_0}, j\neq j_0\}$. We can choose any system of generators of the diagonal ideal
 $\{h_j\}$,  $j\in J'$, which will define a coordinate system $\{h_j, z_{j_o}\}$ for any $j_0$, and obtain
  \begin{equation}
 \nonumber
 \Delta^{(J)}_\ast(M)\stackrel{\sim}{\rightarrow}M[h_j^{-1},\; j\in J']\bigwedge_{j\in J'}\frac{dh_j}{h_j}.
 \end{equation}
 The function that an element $m\otimes \bigwedge_{j\in J'}dh_j/h_j^{-k_j-1}$ represents is as follows: restrict a function $f$ to the formal neighborhood of the diagonal, write it as a power series in $h_j$, $j\in J'$, with coefficients in $\BC[z_{j_0}]$, then  $m\otimes \bigwedge_{j\in J'}dh_j/h_j^{-k_j-1}$ will pick the coefficient of  $\prod_{j\in J'}h_j^{k_j}$ applied to $m$.

 \subsection{  }
 \label{d-mod-gener-right}
 
In the case of left $D_C$-modules, we  will be interested not so much in a similar standard direct image functor  as in the functor 
$\widehat{\Delta}^{(J/I)}_*$  defined in \cite{BD}, sect. 3.5,
 to be the right adjoint to the usual inverse image functor $\Delta^{(J/I)*}$. One has
 \begin{equation}
 \label{constrofbddirimage}
 \Gamma(U, \widehat{\Delta}^{(J/I)}_*M)=\mbox{Hom}(\Delta^{(J/I)*}D_U,M).
 \end{equation}
 A receptacle for an OPE map is defined by localization:
 \[
 \widetilde{\Delta}^{(J/I)}_*M=j_*^{(J/I)}\CO_{U^{(J/I)}}\otimes_{\CO_{C^J}}\widehat{\Delta}^{(J/I)}_*M,
 \]
 where $j^{(J/I)}:U^{(J/I)}\rightarrow C^J$ is the complement to diagonals that contain $\Delta^{(J/I)}(C^I)\subset C^J$.
 
 If $C$ is the line $\BC$, $I=\{1,2,3,...,n\}$, then $M$ can be replaced with the space of global sections and (\ref{constrofbddirimage}) suggests a realization
 \[
 \widehat{\Delta}^{(I)}_*M\stackrel{\sim}{\rightarrow} pr^*_{i_0}M[[(z_i-z_{i_0}),\;i\neq i_0]],
 \]
 where $pr_{i_0}: C^I\rightarrow C$ is the projection number $i_0$. Of course, there are $n$ such realizations, one for each value of $i_0$, but they are pairwise canonically identified, as $D$-modules,  by a formal application of the Taylor expansion:
 \begin{equation}
 \label{pairwise ident}
 m(z_{i_0})\stackrel{\sim}{\rightarrow} \sum_{n=0}^\infty \frac{(z_{i_0}-z_{i_1})^n}{n!}\partial_{z_{i_1}} ^n m(z_{i_1}),
 \end{equation}
 which is nothing but Grothendieck's way of defining a $D$-module structure.
 Here $m(z_{i_j})$ means not a function of $z_{i_j}$, but $m\in M$ regarded as a section of $pr_{i_j}^*M$, the convention we will be using repeatedly.
 
 By the same token,  we can say that $ \widetilde{\Delta}^{(I)}_*M$ is an $I$-tuple of modules
  \[
 pr^*_{i_0}M[[(z_i-z_{i_0}),\;i\neq i_0]][(z_i-z_j)^{-1},\;i\neq j],\; i_0\in I,
 \]
 along with pairwise identifications defined by  (\ref{pairwise ident}).

\subsection{ } 
\label{lie-star-alg}
Beilinson and Drinfeld \cite{BD} promote $\mbox{Mod}^r(D_C)$ to a  pseudo-tensor category with the following spaces of *-operations: 
\[
P^*_I(\{M_i\},N)=Hom_{D_{C^I}}(\boxtimes_I M_i,\Delta^{(I)}_*N).
\]
Such operations can be composed: given a surjection $\pi: J\rightarrow I$, a $J$-family of objects $\{L_j\}$, an $I$-family of operations
$\phi_i\in P^*_{J_i}(\{L_j\},M_i)$, where $J_i=\pi^{-1}(i)$, and $\psi\in P^*_I(\{M_i\},N)$, the composition $\psi(\{\phi_i\})\in P^*_J(\{L_j\},N)$ is defined to be the following composite map
\[
\boxtimes_JL_j=\boxtimes_I\boxtimes_{J_i}L_j\stackrel{\boxtimes_I\phi_i}{\longrightarrow}\boxtimes_I\Delta^{(J_i)}_*M_i=\Delta^{(J/I)}_*(\boxtimes_IM_i)
\stackrel{\Delta^{(J/I)}_*(\psi)}{\longrightarrow} \Delta^{(J/I)}_*\circ\Delta^{(I)}_*N\stackrel{\sim}{\rightarrow}\Delta^{(J)}_*N.
\]

By definition, \cite{BD}, a Lie* algebra is a Lie algebra object in this pseudo-tensor category; that is to say, a Lie* algebra is  a pair of a $D_C$-module $L$ and an operation $[.,.]\in P^*_2(\{L,L\},L)$, a Lie* bracket, which is anticommutative and satisfies the Jacobi identity. Let us now see what this means.

 \subsection{ }
 \label{expl-form-lie-star}
 We will be exclusively working in the situation where  $C=\BC$ and all the structure elements are translation-invariant. This means that a left $D_\BC$-module is of the form $V\otimes_\BC\CO_\BC$ for some vector space $V$, the $D_\BC$-module structure is determined by a linear map $\partial:V\rightarrow V$ so that $\partial_z(v\otimes f(z))=\partial v\otimes f(z)+v\otimes f(z)'$.  Any right $D_\BC$-module is $V\otimes_{\BC}\Omega_\BC$
 for some left $D_\BC$-module $V$ as  above with action $(v\otimes f(z)\,dz)\partial_z=-\partial v\otimes f(z)\,dz-v\otimes f(z)'\,dz$. We will write $v\partial$ for $-\partial v$.
 
 Similarly, being an $\CO_{C\times C}$-morphism, a binary operation
 \[
 [.,.]:\; V\otimes_{\BC}\Omega_\BC\boxtimes V\otimes_{\BC}\Omega_\BC\rightarrow \Delta^{(\{1,2\})}_* (V\otimes_{\BC}\Omega_\BC)
 \]
 is determined by its restriction to the fiber $Vdz_1\otimes_\BC Vdz_2$, and its translation-invariance means this restriction has the form
 \[
 [adz_1,bdz_2]=\sum_{n=0}^{\infty} a_{(n)}b (z_2)(z_1-z_2)^{-n-1}dz_1\wedge dz_2, \; 
  \]
  where $a_{(n)}b\in V$, $a_{(n)}b=0$ for all $n$ large enough; we are using the functional realization (\ref{functional-realiz}).  Hence $[.,.]$ is determined by an infinite family of  products $_{(n)}:V\otimes_{\BC}V\rightarrow V$. 
  
  Being a $D_{C\times C}$-morphism means the following 2 equalities of elements of $\mbox{End}_\BC(V)$:
  \begin{equation}
  \label{expl form lin over d}
  (\partial a)_{(n)}=[\partial,a_{(n)}]=-na_{(n-1)}.
  \end{equation}
  
  Let us spell out the anticommutativity and the Jacobi identity in this setting.
  
  {\em Anticommutativity:}
  
  If 
   \[
 [a(z_1)dz_1,b(z_2)dz_2]=\sum_{n=0}^{\infty} a_{(n)}b (z_2)(z_1-z_2)^{-n-1}dz_1\wedge dz_2, \; 
  \]
then
 \begin{eqnarray*}
 [b(z_2)dz_2,a(z_1)dz_1]&=&\sum_{n=0}^{\infty} b_{(n)}a (z_1)(z_2-z_1)^{-n-1}dz_2\wedge dz_1\\
 &=&
 \sum_{n=0}^{\infty} \sum_{j=0}^\infty \frac{(z_1-z_2)^j}{j!}\partial^j(b_{(n)}a )(z_2)(z_2-z_1)^{-n-1}dz_2\wedge dz_1,
  \end{eqnarray*}
 (the last equality is the result of an application of (\ref{pairwise ident}).) 
  Comparing the coefficients in front of $(z_1-z_2)^{-n-1}$ we see that
 
  \begin{eqnarray}
   [a(z_1)dz_1,b(z_2)dz_2]&=&-[b(z_2)dz_2,a(z_1)dz_1] \Leftrightarrow\nonumber\\
   a_{(n)}b&=&(-1)^{n+1}\sum_{j=0}^\infty\frac{(-1)^j}{j!}\partial^j(b_{(n+j)}a)\mbox{ for all }n,\label{anticomm-in-n-prod}
  \end{eqnarray}
  the identity long known in vertex algebra theory as {\em skew-commutativity.} 
  
  {\em Jacobi:} 
  
  The definition of the composition of *-operations, sect. \ref{lie-star-alg}, gives
  \begin{eqnarray*}
 & & [a(z_1)dz_1,[b(z_2)dz_2,c(z_3)dz_3]]=[a(z_1)dz_1,\sum_{m=0}^\infty b_{(m)}c(z_3)(z_2-z_3)^{-m-1}dz_2\wedge dz_3]\\
  &=&\sum_{m,n=0}^\infty a_{(n)}(b_{(m)}c)(z_3)(z_1-z_3)^{-n-1}(z_2-z_3)^{-m-1}dz_1\wedge dz_2\wedge dz_3;
  \end{eqnarray*}
  similarly,
  \begin{equation*}
   [b(z_2)dz_2,[a(z_1)dz_1,c(z_3)dz_3]]
  =\sum_{m,n=0}^\infty b_{(m)}(a_{(n)}c)(z_3)(z_1-z_3)^{-n-1}(z_2-z_3)^{-m-1}dz_1\wedge dz_2\wedge dz_3;
  \end{equation*}
  and
  \begin{equation*}
 [[a(z_1)dz_1,b(z_2)dz_2],c(z_3)dz_3]=\sum_{s,t=0}^\infty (a_{(s)}b)_{(t)}c(z_3) (z_1-z_2)^{-s-1}(z_2-z_3)^{-t-1}dz_1\wedge dz_2\wedge dz_3.
 \end{equation*}
 
 We will show that the Jacobi identity,
 \[
 [a(z_1)dz_1,[b(z_2)dz_2,c(z_3)dz_3]]- [b(z_2)dz_2,[a(z_1)dz_1,c(z_3)dz_3]]= [[a(z_1)dz_1,b(z_2)dz_2],c(z_3)dz_3],
 \]
 amounts to  the well-known {\em Borcherds commutator formula}:
 \begin{equation}
 \label{borch-comm-in-n-ope}
 [a_{(n)},b_{(m)}]=\sum_{s=0}^\infty{n\choose s}(a_{(s)}b)_{(n+m-s)};
 \end{equation}
 in fact, we will do this in two different ways.
 
 Guided by sect. \ref{d-mod-gener}, we need to show that (\ref{borch-comm-in-n-ope}) is equivalent to the corresponding differential operators being equal to each other. {\em One way} to do this is to choose an arbitrary test function, say $(z_1-z_3)^n(z_2-z_3)^m$, expand it in powers of the chosen coordinate system, $\{(z_1-z_3), (z_2-z_3), z_3\}$ for the first two expression and $\{(z_1-z_2), (z_2-z_3),z_3\}$ for the third one, and then take the residues.  We obtain:
 
 $[a(z_1)dz_1,[b(z_2)dz_2,c(z_3)dz_3]]$ sends  $(z_1-z_3)^n(z_2-z_3)^m$ to $a_{(n)}(b_{(m)}c)$;
 
 $[b(z_2)dz_2,[a(z_1)dz_1,c(z_3)dz_3]]$ sends  $(z_1-z_3)^n(z_2-z_3)^m$ to $b_{(m)}(a_{(n)}c)$; and 
 
  $[[a(z_1)dz_1,b(z_2)dz_2],c(z_3)dz_3]$  sends  $(z_1-z_3)^n(z_2-z_3)^m$  to $\sum_{s=0}^\infty{n\choose s}(a_{(s)}b)_{(n+m-s)}$; this is  because we have to re-expand as follows:
  \[
  (z_1-z_3)^n(z_2-z_3)^m=\sum_{s=0}^\infty{n\choose s}(z_1-z_2)^s(z_2-z_3)^{n+m-s}.
  \]
  That the 1st value minus the 2nd equals the 3rd for all $m,n$ is the content of (\ref{borch-comm-in-n-ope}).
  
  {\em Allternatively,} and more straightforwardly, let us compare the differential operators themselves.  As explained in the end of sect. \ref{d-mod-gener}, this amounts to choosing the coordinate systems, the same we have just chosen, and the corresponding dual bases of the tangent sheaf; these are $\{\partial_{z_1},\partial_{z_2},
  \partial_{z_1}+ \partial_{z_2}+ \partial_{z_3}\}$ for the first two expressions and  $\{\partial_{z_1},\partial_{z_1}+\partial_{z_2},
  \partial_{z_1}+ \partial_{z_2}+ \partial_{z_3}\}$ for the third.

 We obtain
  \begin{equation*}
 [a(z_1)dz_1,[b(z_2)dz_2,c(z_3)dz_3]]
  =\sum_{m,n=0}^\infty a_{(n)}(b_{(m)}c)(z_3)dz_3\otimes \frac{\partial_{z_1}^n}{n!}\frac{\partial_{z_2}^{m}}{m!},
  \end{equation*}
   \begin{equation*}
 [b(z_2)dz_2,[a(z_1)dz_1,c(z_3)dz_3]]
  =\sum_{m,n=0}^\infty b_{(m)}(a_{(n)}c)(z_3)dz_3\otimes \frac{\partial_{z_1}^n}{n!}\frac{\partial_{z_2}^{m}}{m!},
  \end{equation*}
  and
  \begin{equation*}
 [[a(z_1)dz_1,b(z_2)dz_2],c(z_3)dz_3]=\sum_{s,t=0}^\infty (a_{(s)}b)_{(t)}c(z_3)dz_3\otimes \frac{\partial_{z_1}^{s}}{s!}\frac{(\partial_{z_1}+\partial_{z_2})^{t}}{t!}.
 \end{equation*}
 
 Comparing coefficients in front of $\partial_{z_1}^n\partial_{z_2}^m/n!m!$ we again obtain (\ref{borch-comm-in-n-ope}). $\qed$
 
 \subsubsection{Remark.}
 \label{a third proof by expansion}
 A curious reader might enjoy seeing yet another, a third,  proof of (\ref{borch-comm-in-n-ope}). One needs to take the third term of the Jacobi identity,
  \begin{equation*}
 [[a(z_1)dz_1,b(z_2)dz_2],c(z_3)dz_3]=\sum_{s,t=0}^\infty (a_{(s)}b)_{(t)}c(z_3) (z_1-z_2)^{-s-1}(z_2-z_3)^{-t-1}dz_1\wedge dz_2\wedge dz_3,
 \end{equation*}
 and re-expand in powers of $(z_1-z_3)$ and $(z_2-z_3)$  making sure that the term $(z_1-z_2)^{-s-1}$ is expanded in positive powers of $(z_2-z_3)$. The result is this:
 \begin{eqnarray*}
 (z_1-z_2)^{-s-1}(z_2-z_3)^{-t-1}&=&(z_1-z_3)^{-s-1}(1-\frac{z_2-z_3}{z_1-z_3})^{-s-1}(z_2-z_3)^{-t-1}\\
 &=& \sum_{j=0}^\infty {j\choose s}(z_1-z_3)^{-j-1}(z_2-z_3)^{-t+j-s-1}.
 \end{eqnarray*}
 The familiar device of comparing the coefficients in front of  $(z_1-z_3)^{-n-1}(z_2-z_3)^{-m-1}$ (and ignoring the terms with positive powers) will give the desired (\ref{borch-comm-in-n-ope}).
 
 The reason this works is as follows.  Consider $C\times C$ along with the diagonal $\Delta^{(2)}:C\rightarrow C\times C$ and its complement
 $j: U\rightarrow C\times C$. Then there is a residue map
 \[
 res: j_*j^* M_1\boxtimes M_2\rightarrow \Delta^{(2)}_*(M_1\otimes^!M_2),
 \]
 where $M_1\otimes^!M_2=M_1\otimes_{\CO_C}(M_2\otimes_{\CO_C}\Omega_C^{-1})$, see \cite{BD}, especially sect. 3.2. Computationally, this maps amount to choosing local coordinates, say $z_1-z_2, z_2$, expanding in powers of $z_1-z_2$ with coefficients depending on $z_2$, and keeping only singular, w.r.t. $z_1-z_2$, terms.
 
 This construction has a higher dimensional analogue. Given a finite set $T$, consider the complement to all diagonals, $j^{(T)}: U^{(T)}\rightarrow C^T$.
 Given a surjection $\pi: T\rightarrow T'$ with $|T|=|T'|+1$, there naturally and similarly arises a residue map
  \[
 res_{T/T'}: j^{(T)}_*j^{(T)*} \boxtimes_T M_t\rightarrow \Delta^{(T/T')}_*j^{(T')}_*j^{(T')*}\boxtimes_{T'}\otimes^!_{\pi^{-1}(t')}M_t.
 \]
 This construction can be  also iterated: given 2 surjections $T\stackrel{\pi_1}{\rightarrow} T'\stackrel{\pi_2}{\rightarrow} T''$, $|T|=|T'|+1=|T''+2$, one can consider the composition
  \begin{eqnarray*}
 j^{(T)}_*j^{(T)*} \boxtimes_T M_t&\stackrel{res_{T/T'} }{\longrightarrow}& \Delta^{(T/T')}_*j^{(T')}_*j^{(T')*}\boxtimes_{T'}\otimes^!_{\pi_1^{-1}(t')}M_t\\
 &\stackrel{\Delta^{(T/T')}(res_{T'/T''} )}{\longrightarrow} &\Delta^{(T/T')}_*( \Delta^{(T'/T'')}_*j^{(T'')}_*j^{(T'')*}\boxtimes_{T''}\otimes^!_{(\pi_2\circ\pi_1)^{-1}(t')}M_t)\\
 &\stackrel{\sim}{\longrightarrow}&\Delta^{(T/T'')}_*j^{(T'')}_*j^{(T'')*}\boxtimes_{T''}\otimes^!_{(\pi_2\circ\pi_1)^{-1}(t')}M_t
 \end{eqnarray*}
 
 The suggested proof of  (\ref{borch-comm-in-n-ope}) deals with such a composition where $T=\{1,2,3\}$, $T'=\{1,\{2,3\}\}$, which means a 2-element partition of $T$), and $T''=\{\{1,2,3\}\}$, a 1-element set;  the maps $T\rightarrow T'\rightarrow T''$ are the  obvious  ones; and $M_1=M_2=\Omega_C$, $M_3=V\otimes\Omega_\BC$. Under this composition
 \begin{eqnarray*}
  a_{(n)}(b_{(m)}c)(z_3)(z_1-z_3)^{-n-1}(z_2-z_3)^{-m-1}dz_1\boxtimes dz_2\boxtimes dz_3&\mapsto&\\
   a_{(n)}(b_{(m)}c)(z_3)(z_1-z_3)^{-n-1}(z_2-z_3)^{-m-1}dz_1\wedge dz_2\wedge dz_3,& &\\
   b_{(m)}(a_{(n)}c)(z_3)(z_1-z_3)^{-n-1}(z_2-z_3)^{-m-1}dz_1\boxtimes dz_2\boxtimes dz_3&\mapsto&\\
   b_{(m)}(a_{(n)}c)(z_3)(z_1-z_3)^{-n-1}(z_2-z_3)^{-m-1}dz_1\wedge dz_2\wedge dz_3,& &\\
   (a_{(s)}b)_{(t)}c(z_3) (z_1-z_2)^{-s-1}(z_2-z_3)^{-t-1}dz_1\boxtimes dz_2\boxtimes dz_3&\mapsto&\\
   (a_{(s)}b)_{(t)}c(z_3) (z_1-z_2)^{-s-1}(z_2-z_3)^{-t-1}dz_1\wedge dz_2\wedge dz_3,
   \end{eqnarray*}
   a somewhat uninspiring result, but notice that although the images appear indistinguishable from the preimages, they belong to different spaces and are of a different nature; e.g. the elements on the right are zero divisors, those on the left are rational functions. The first two are obtained using $z_1-z_3,z_2-z_3,z_3$ as a coordinate system, the last one uses $z_1-z_2,z_2-z_3,z_3$. The computation carried out at the beginning of this remark simply rewrites the last one using $z_1-z_3,z_2-z_3,z_3$.

 \subsection{ }
 \label{compound pseudotensor str}
 The usual tensor product of left $D$-modules, $M\otimes_{\CO_C}N$, has an analogue for the right $D$-modules, \cite{BD}:
 \[
 M\otimes^!N=M\otimes_{\CO_C}(N\otimes_{\CO_C}\Omega_C^{-1}),
 \]
 which in the translation-invariant setting is nothing but the tensor product of fibers:
  \[
 M[z]dz\otimes^!N[z]dz=(M\otimes_\BC N)[z]dz.
 \]
 Therefore we can consider associative (commutative) algebras relative to the functor $\otimes^!$; in the translation-invariant setting, $A[z]dz$ being such an algebra means simply $A$ being an ordinary associative (commutative) $\BC$-algebra.
 
 Suppose now $P$ carries a Lie* algebra bracket and an associative commutative $\otimes^!$-multiplication. In order to see whether or not these two operations are compatible, one needs an extra structure, called a {\em compound tensor product} of *-operations, \cite{BD}. 
 
 Two surjections,
 $\pi^{(I/S)}:I\rightarrow S$ and $\pi^{(I/T)}:I\rightarrow T$ are called {\em complementary} if $\Delta^{(I/S)}(C^S)$ and $\Delta^{(I/T)}(C^T)$ transversally intersect at the main diagonal, $ \Delta^{(I)}(C)$.  If so, and given
 
 (1)  an $I$ family of right $D$-modules $\{L_i\}$,

 (2)  an $S$-family of right $D$-modules $\{M_s\}$,  and
 
 (3) an $S$-family of *-operations
 \begin{equation*}
 \phi_s: \boxtimes_{I_s}L_i\longrightarrow \Delta^{(I_s)}_*M_s,\; s\in S,
 \end{equation*}
 where $I_s\subset I$ is the preimage of $s\in S$ w.r.t. $\pi^{(I/S)}:I\rightarrow S$,
 
 we define an element of $P^*_T(\{\otimes^!_{I_t}L_i\},\otimes_S^!M_s)$. This element is called a {\em tensor product of operations} and denoted 
 by     $\otimes^{I}_{S,T}(\phi_s)$. In other words, we will define a map

 \[
 \otimes^{I}_{S,T}:\; \otimes_SP^*_{I_s}(\{L_i\},M_s)\rightarrow P^*_T(\{\otimes^!_{I_t}L_i\},\otimes_S^!M_s).
 \]
 This map is defined essentially by pull-back onto $\Delta^{(I/T)}(C^T)$ as follows:
 
 since we have
 \[
 \boxtimes_IL_i=\boxtimes_S\boxtimes_{I_s}\stackrel{\boxtimes_S\phi_s}{\longrightarrow}\boxtimes_S\Delta^{(I_s)}_*M_s=\Delta^{(I/S)}_*\boxtimes_S M_s,
 \]
 we apply the inverse image functor $\Delta^{(I/T)!}$ (properly understood, \cite{BD}, 2.2.6) and define  
 \[
 \otimes^{I}_{S,T}(\phi_s):\; \boxtimes_T(\otimes_{I_t}^!L_i)\longrightarrow \Delta_*^{(T)}(\otimes^!_S M_s)
 \]
  to be the following composition
  \begin{eqnarray*}
\boxtimes_T(\otimes_{I_t}^!L_i)&\stackrel{\sim}{\rightarrow}& \Delta^{(I/T)!}(\boxtimes_IL_i)\stackrel{ \Delta^{(I/T)!}(\boxtimes_S\phi_s)}\longrightarrow
 \Delta^{(I/T)!}\circ\Delta^{(I/S)}_*\boxtimes_S M_s\\
&\longrightarrow& \Delta_*^{(T)}\circ\Delta^{(S)!}(\boxtimes_S M_s)
\stackrel{\sim}{\longrightarrow}\Delta_*^{(T)}(\otimes^!_S M_s),
 \end{eqnarray*}
 the penultimate arrow being the  base change isomorphism, the point where the complementarity is used.
 
 The most important -- and simplest -- case of this construction arises when $S$ is a 2-element set, say $S=\{a,b\}$; it defines a partition
 $I=I_a\sqcup I_b$; having a complementary $I\rightarrow T$ means choosing $i_a\in I_a$ and $i_b\in I_b$ and letting 
 \[
 T=I_a\vee_{(i_a,i_b)}I_b,
 \]
 that is, $I$  with $i_a$ identified with  $i_b$.
  What the above construction
 accomplishes is, given two *-operations, $\phi_{I_a}\in P^*_{I_a}(\{L_i\},M_a)$ and $\phi_{I_b}\in P^*_{I_b}(\{L_i\},M_b)$ it defines a tensor product,
 $\phi_{I_a}\otimes_{i_a,i_b}\phi_{I_b}\in P^*_{I_a\vee_{(i_a,i_b)}I_b}(\{L_i, L_{i_a}\otimes^! L_{i_b},i\neq i_a, i_b\},M_a\otimes^!M_b)$. This tensor product is obtained by pulling the morphism
 \[
 \phi_{I_a}\boxtimes\phi_{I_b}:(\boxtimes_{I_a}L_i)\boxtimes(\boxtimes_{I_b}L_i)\rightarrow \Delta^{(I_a)}_*M_a\boxtimes  \Delta^{(I_b)}_*M_b
 \]
 back onto $C^{I_a\vee_{(i_a,i_b)}I_b}=\{z_{i_a}=z_{i_b}\}\subset C^{I_a\sqcup I_b}$.
 
 Therefore, if $L$ is a Lie* algebra with bracket $[.,.]$,  then we can define its action on $L\otimes^!L$ by derivations: 
 \begin{equation}
 \label{action on tens prod}
 [a,b\otimes^!c]=[a,b]\otimes^!c+b\otimes^![a,c],
 \end{equation}
 where $[a,b]\otimes^!c= ([.,.]_{\{1,2\}}\otimes_{2,3}\mbox{Id}_{\{3\}})(a,b,c)$ and $b\otimes^![a,c]=  \mbox{Id}_{\{2\}}\otimes_{2,3}[.,.]_{\{1,3\}} (a,b,c)$.
 
 Explicitly, and working in the translation-invariant situation, if $[a(z_1)dz_1,b(z_2)dz_2]=\sum_{n\geq 0} a_{(n)}b(z_2)dz_2(z_1-z_2)^{-n-1}$, then
 \begin{equation}
 \label{leib-1}
 [a(z_1)dz_1,b(z_2)dz_2]\otimes^!c(z_2)dz_2]=\sum_{n\geq 0} ((a_{(n)}b)\otimes c)(z_2)dz_1\wedge dz_2(z_1-z_2)^{-n-1},
 \end{equation}
 \begin{equation}
 \label{leib-2}
b(z_2)dz_2\otimes^! [a(z_1)dz_1,c(z_2)dz_2]=\sum_{n\geq 0} (b\otimes(a_{(n)}c))(z_2) dz_1\wedge dz_2(z_1-z_2)^{-n-1}.
 \end{equation}
 It is curious to note that thanks to the choices we have made changing the order of the arguments leads to a formula, which looks drastically different from
 (\ref{leib-1}, \ref{leib-2}):
 \begin{eqnarray}
  [b(z_2) dz_2,a(z_1)dz_1]\otimes^!c(z_2)&=&\sum_{n\geq 0} (b_{(n)}a)(z_1)\otimes c(z_2)dz_1\wedge dz_2(z_2-z_1)^{-n-1}\nonumber\\
  &=& \sum_{n\geq j} ((b_{(n)}a)\otimes \frac{1}{j!}\partial^jc)(z_1)dz_1\wedge dz_2(z_2-z_1)^{-n-1+j},\label{leib-3}
  \end{eqnarray}
  and similarly
   \begin{equation}
   \label{leib-4}
  b(z_2) dz_2\otimes^![c(z_2),a(z_1)dz_1]
  = \sum_{n\geq j} \frac{1}{j!}(\partial^jb\otimes (c_{(n)}a))(z_1)dz_1\wedge dz_2(z_2-z_1)^{-n-1+j},
  \end{equation}
  
  \subsection{ Definition.}
  \label{def of coiss alg}
 {\em A right $D_C$-module $A$ is a coisson algebra if it is a Lie* algebra with bracket $\{.,.\}$ and a unital commutative associative $\otimes^!$-algebra such that the Leibniz rule holds, that is, the product
  \[
  A\otimes^!A\rightarrow A
  \]
  is a morphism of $A$-modules, where on the right $A$ is regarded as the adjoint module and on the left the $A$-module structure is defined by 
  (\ref{action on tens prod}).}
  
  Of course, in the translation-invariant setting this clear-cut definition can be made quite explicit but opaque, whereby a vertex Poisson algebra is a unital, associative $\BC$-algebra $A$ with derivation $\partial$ and a family of products $_{(n)}:A\otimes A\rightarrow A$, $n\geq 0$, such that the latter satisfy the identities (\ref{expl form lin over d},\ref{anticomm-in-n-prod},\ref{borch-comm-in-n-ope}) and  are derivations of $A$ as an associative algebra: 
  \begin{equation}
  \label{leib-6}
  a_{(n)}(bc)=(a_{(n)}b)c+b(a_{(n)}c). 
  \end{equation}
  
  The latter identity uses an obvious restatement of (\ref{leib-1},\ref{leib-2}). Incidentally, a similar restatement  of (\ref{leib-3},\ref{leib-4}) gives another form of the Leibniz identity:
  \begin{equation}
  \label{leib-5}
  (bc)_{(n)}a=\sum_{j=0}^\infty \frac{1}{j!}(\partial^jb)(c_{(n+j)}a)+\frac{1}{j!}(\partial^jc)(b_{(n+j)}a).
  \end{equation}
  It  will be seen to be  a quasiclassical limit of the  normal ordering formula (\ref{norm-orde-formu}).

  As explained in the inspirationally written introduction to \cite{BD} the whole of ``chiral theory'' might owe its existence to the following simple observation:
  if $A$ is a Poisson algebra, then the jet algebra $J_\infty A$ is not naturally a Poisson algebra; instead $J_\infty A\otimes \Omega_\BC$  is naturally a coisson  (a.k.a. vertex Poisson ) algebra.  Namely, $J_\infty A$ contains $A$, carries a derivation $\partial: J_\infty A\rightarrow J_\infty A$, and is generated as an algebra by
   the derivatives $\partial^ja$, $a\in A$. If $\{.,.\}$ is A Poisson bracket on $A$, then we define the Lie* bracket by letting
   \[
   \{a(z_1)dz_1,b(z_2)dz_2\}^*=\{a,b\}(z_2)\frac{dz_1\wedge dz_2}{z_1-z_2},\; a,b\in A\subset J_\infty A,
   \]
   \[
   \{(a(z_1)dz_1)\partial ^n ,(b(z_2)dz_2)\partial^m\}^*=
   (\{a,b\}(z_2)\frac{dz_1\wedge dz_2}{z_1-z_2})\partial_{z_1}^n\partial_{z_2}^m,\; a,b\in A\subset J_\infty A,
   \]
   and then extending to the whole of $J_\infty A\otimes \Omega_\BC$ by using the Leibniz identity.
   
    \subsection{Example.}
   \label{Example q class}

  We return to the setup of sect. \ref{Calc1.} and 
   let $p: X\rightarrow S$ be a smooth morphism of smooth algebraic varieties over $\BC$. The relative cotangent bundle $T^*X/S$ is a Poisson variety and
   $\CO_{T^*X/S}$ is a sheaf of Poisson algebras. By passing to the  jet spaces we obtain a sheaf of coisson algebras $\CO_{J_\infty T^*X/S}\otimes\Omega_\BC$.
   
   Everything we have discussed has a super-version. Thus we obtain a Poisson superalgebra $\CO_{T^*\Pi TX/S}$ and its coisson version
   $\CO_{J_\infty T^*\Pi TX/S}\otimes\Omega_\BC$; here $\Pi TX/S$ means $TX/S$ with  fibers of the projection $TX/S\rightarrow X/S$ declared odd.
   
   Working locally allows us to write down explicit formulas. Locally means $X/S$ is affine and
     possesses a  coordinate system, i.e., functions $x_1,...,x_N$ and
   an abelian basis of $\CT_{X/S}$, $\partial_{x_1},\partial_{x_2},...,\partial_{x_N}$ so that $\partial_{x_i}(x_j)=\delta_{ij}$.  This implies that
   $J_\infty\BC[T^*X/S]=\BC[T^*X/S][\xni,\partxni;1\leq i\leq N,n>0] $.
   
   In the supercase another pair of families of, this time odd, variables is adjoined, and we obtain:
   $J_\infty\BC[T^*\Pi TX/S]=\BC[T^*\Pi TX/S][\xni,\partxni, \dxni,\partdxni;1\leq i\leq N,n>0] $.
   
   Both these algebras carry a canonical derivation $\partial$, and this is reflected in the appearance of the upper index: we let
   \begin{equation}
   \label{action of d on jet}
   a^{(n)}=\frac{1}{n!}\partial^n(a)
   \end{equation}
    for $a=x_i, \partial_{x_i}, dx_i, \partial_{dx_i}$, $n\geq 0$.
   
  According to sect. \ref{def of coiss alg}, the ordinary Poisson (super)brackets, $\{\partial_{x_i},x_j\}=\delta_{ij}$ and $\{\partial_{dx_i},dx_j\}=\delta_{ij}$
define the coisson bracket   on $J_\infty\BC[T^*\Pi TX/S]\otimes\Omega_\BC$ as follows:
\begin{eqnarray*}
\{\partxni(z_1) dz_1,x^{(m)}_j(z_2) dz_2\}&=&(-1)^n {n+m\choose n}\delta_{ij}\frac{dz_1\wedge dz_2}{(z_1-z_2)^{n+m+1}},\\
\{\partdxni(z_1) dz_1,dx^{(m)}_j(z_2) dz_2\}&=&(-1)^n {n+m\choose n}\delta_{ij}\frac{dz_1\wedge dz_2}{(z_1-z_2)^{n+m+1}},
\end{eqnarray*}
and more generally, as a repeated application of the Leibniz identity (\ref{leib-6}), (\ref{leib-5}) will show, if $P,Q\in J_\infty\BC[T^* X/S]$
\begin{eqnarray*}
\{P(z_1)dz_1, Q(z_2)dz_2\} = 
\end{eqnarray*}
\begin{eqnarray*}
 =  
\sum_{i,n,m}\sum_{j=0}^\infty\frac{(z_1-z_2)^j}{j!}\partial^j(\frac{\partial P}{\partial(\partxni)})(z_2)\frac{\partial Q}{\partial(x^{(m)}_i)}\{\partxni(z_1),x^{(m)}_i(z_2)\}
\end{eqnarray*}
\begin{eqnarray*}
+ \sum_{i,n,m}\sum_{j=0}^\infty\frac{(z_1-z_2)^j}{j!}\partial^j(\frac{\partial P}{\partial(x^{(m)}_i)})(z_2)\frac{\partial Q}{\partial(\partxni )}
\{x^{(m)}_i(z_1),\partxni(z_2)\}
\end{eqnarray*}

In the superversion  of the above formula  
two more similar summands, those related to $dx_i$ and $\partial_{dx_i}$, will appear, but we will refrain from  writing them down as an attempt to fix signs is only likely to sow confusion. What is important, though, is the fact that such brackets are defined by means of certain universal bi-differential operators with constant coefficients acting on the (commutative, associative) algebra $J_\infty\BC[T^*\Pi TX/S]$; this will prove essential for the quantization of this bracket, the problem to which we now turn.

\subsection{ }
\label{jets in lie alg current alg}
Similarly and more simply, given a Lie algebra $\fg$ we consider 
$$
J_\infty\fg=\fg\otimes_\BC\BC[\partial],
$$
a free $\BC[\partial]$-module generated by $\fg$ and endow  $J_\infty\fg\otimes\Omega_\BC$ with a Lie* algebra structure by letting

 \[
   \{a(z_1)dz_1,b(z_2)dz_2\}^*=[a,b](z_2)\frac{dz_1\wedge dz_2}{z_1-z_2},\; a,b\in \fg=\fg\otimes 1\subset J_\infty \fg,
   \]
   \[
   \{(a(z_1)dz_1)\partial ^n ,(b(z_2)dz_2)\partial^m\}^*=
   ([a,b](z_2)\frac{dz_1\wedge dz_2}{z_1-z_2})\partial_{z_1}^n\partial_{z_2}^m,\; a,b\in \fg\subset J_\infty \fg,
   \]
   For example, a moment's thought will show that  $J_\infty\BC[T^*\Pi TX/S]\otimes\Omega_\BC$ regarded as a Lie* algebra contains a Lie* subalgebra $J_\infty \CT_{X/S}\otimes\Omega_\BC$.

\subsection{ }
\label{chiral algebra}
The concept of a coisson algebra $V$ is a quasiclassical limit of that of a chiral algebra, a.k.a. a vertex algebra in the translation-invariant situation. By definition, \cite{BD},
a chiral algebra is a Lie algebra object of $\mbox{Mod}^r(D_C)$ with a pseudo-tensor structure  defined by means of the spaces of chiral operations,
\[
P^{ch}_I(\{M_i\},N)\stackrel{\mbox{def}}{=}\mbox{Hom}_{D_{C^I}}(j_*^{(I)}j^{(I)*}(\boxtimes_IM_i),\Delta^{(I)}_*N),
\]
where $j^{(I)}:U^{(I)}\hookrightarrow C^I$ is the embedding of the open subset $U^{(I)}=\{(z_i)\mbox{ s.t. } z_{i_1}\neq z_{i_2}\mbox{ if }i_1\neq i_2\}$, cf.   sect. \ref{lie-star-alg}; it also contains a unit, $\Omega_C\rightarrow V$, whose definition we will skip.

A chiral bracket  on $M\in \mbox{Mod}^r(D_C)$ is then a mophism
\[
[.,.]:j_*^{(2)}j^{(2)*} M\boxtimes M\rightarrow \Delta^{(2)}_*M,
\]
which is antisymmetric and satisfies the Jacobi identity.

Computationally more convenient than a chiral bracket is the concept of an operator product expansion, OPE for short. It is easy to understand, \cite{BD}, 3.5, that
a morphism $j_*^{(2)}j^{(2)*} M\boxtimes M\rightarrow \Delta^{(2)}_*M$
  is the same thing as a morphism of left $D_{C^2}$-modules
\[
\circ: M^l\boxtimes M^l\rightarrow \widetilde{\Delta}^{(2)}_*M^l,
\]
where $M^l=M\otimes_{\CO_C}\Omega_C^{-1}$ and  $ \widetilde{\Delta}^{(2)}_*M^l$ is the $D_{C^2}$-module that was defined in sect. \ref {d-mod-gener}.

In the translation-invariant setting, where $M=V\otimes_\BC\Omega_{\BC}$, $M^l=V\otimes_\BC\CO_{\BC}$, the OPE is determined by its restriction to the fiber:
\begin{equation}
\label{what ope looks like in tr inva}
 pr_1^*V\otimes _\BC pr_2^*V\ni a\otimes b\mapsto a(z_1)\circ b(z_2)=\sum_{n=-\infty}^\infty a_{(n)}b(z_2)(z_1-z_2)^{-n-1}\in V((z_1-z_2)),
 \end{equation}
 and the $\BC[\partial_{z_1},\partial_{z_2}]$-linearity is equivalent to the familiar identity (\ref{expl form lin over d}), except now it is valid for all $n$, positive and negative:
  \begin{equation}
  \label{expl form lin over d all n}
  (\partial a)_{(n)}=[\partial,a_{(n)}]=-na_{(n-1)},\; n\in\BZ.
  \end{equation}
 
 The corresponding chiral   bracket   is as follows
 \begin{eqnarray}
 \label{chiral bracket formula tr inva}
 \frac{adz_1\otimes b dz_2}{(z_1-z_2)^M}&\mapsto& \sum_{n=-M}^\infty a_{(n)}b(z_2)(z_1-z_2)^{-n-1-M} dz_1\wedge dz_2\\
 &\in& V((z_1-z_2)) dz_1\wedge dz_2 \mbox{ mod } V[[z_1-z_2]] dz_1\wedge dz_2.\nonumber,
 \end{eqnarray}
 which, by the way is an illustration of how the OPE point of view is computationally convenient.
 
 The anticommutativity and the Jacobi identity of a chiral bracket are equivalent to the commutativity and associativity of the OPE $\circ$, \cite{BD}, 3.5.10.

 In terms of products $_{(n)}$, commutativity of $\circ$, equivalently anticommutativity of $[.,.]$, means the familiar identity (\ref{anticomm-in-n-prod}), except now it is valid for all $n\in\BZ$, not necessarily positive.  The Jacobi identity for $[.,.]$ is not easy to write down explicitly because of the poles, harder at least  than the computation in sect. \ref{expl-form-lie-star}. Even defining the concept of  associativity for $\circ$ requires an effort. We have
 \begin{eqnarray*}
 a(z_1)\circ(b(z_2)\circ c(z_3))&=&\sum_{m,n} a_{(n)}(b_{(m)}c)(z_3)(z_1-z_3)^{-n-1}(z_2-z_3)^{-m-1}\\
 &\in& V((z_2-z_3))((z_1-z_3)),
 \end{eqnarray*}
  \begin{eqnarray*}
( a(z_1)\circ b(z_2))\circ c(z_3)&=&\sum_{s,t}(( a_{(s)}b)_{(t)}c)(z_3)(z_1-z_2)^{-s-1}(z_2-z_3)^{-t-1}\\
&\in& V((z_1-z_2))((z_2-z_3)),
 \end{eqnarray*}
 two expansions which belong to different spaces. The OPE $\circ$ is called {\em associative} if they are the Laurent series expansions, in their respective domains, of one and the same element of the localization $V[[z_1,z_2,z_3]][(z_1-z_2)^{-1}, (z_1-z_3)^{-1}, (z_2-z_3)^{-1}]$. For example, the associativity demands that in order to obtain $ a(z_1)\circ(b(z_2)\circ c(z_3))$, one has to expand $1/(z_1-z_2)$ as this:
 \[
 \frac{1}{z_1-z_2}=\sum_{n=0}^\infty\frac{(z_2-z_3)^n}{(z_1-z_3)^{n+1}}.
 \]
 None of this is very explicit. One obtains an explicit formula -- the celebrated Borcherds identity -- if one passes to the de Rham cohomology or, equivalently, integrates. This can be done even at the level of a Lie* algebra, the point where we start; we will return to the chiral algebra case in sect.~\ref{subs on borch ident2}.
 
  \subsection{ }
 \label{subs on borch ident1}
 If $M$ is a right $D_C$-module, then one considers its de Rham cohomology $\CH^1_{dR}(C,M)=M/M\CT_C$.
 
 Given a morphism $M\boxtimes M\rightarrow \Delta^{(2)}_*M$, by taking the  de Rham cohomology along the 1st  copy of $C$, one obtains a map
 \[
 \CH^1_{dR}(C,M)\otimes_\BC M\rightarrow M,
 \]
 and by taking the de Rham cohomology over $C\times C$ one obtains a map
 \[
 \CH^1_{dR}(C,M)\otimes_\BC \CH^1_{dR}(C,M)\rightarrow \CH^1_{dR}(C,M).
 \]
 If $M$ is a Lie* algebra (in particular, a chiral algebra), then the latter map is a Lie algebra bracket and the former an action of (the sheaf of Lie algebras)
  $\CH^1_{dR}(C,M)$ on $M$.

 For example, in the translation-invariant setting the answer for the space of sections of $\CH^1_{dR}(C,V((z))dz)$ over a formal disc is as follows  
 \begin{equation*}
 \mbox{Lie}(V)=V((z))dz/\langle (mz^{m-1}a+z^m\partial a)dz,\; a\in V, m\in \BZ\rangle
 \end{equation*}
 
 The chiral bracket formula (\ref{chiral bracket formula tr inva}) shows that the action map is as follows
 \begin{eqnarray}
 \label{acti of lie v - tr inv}
  \mbox{Lie}(V)\otimes V((z))dz&\rightarrow& V((z))dz\\
  \overline{f(z_1) a(z_1) dz_1}\otimes g(z_2) b(z_2) dz_2&\mapsto&\sum_{j=0}^\infty \frac{1}{j!}f^{(j)}(z_2)g(z_2)(a_{(j)}b)(z_2) dz_2\nonumber.
  \end{eqnarray}
 (Indeed, by construction we have to write
 \[
f(z_1)a(z_1)dz_1\otimes g(z_2)b(z_2) dz_2\mapsto \sum_{j=0}^\infty a_{(j)}b(z_2)\frac{f(z_1)g(z_2)}{(z_1-z_2)^{j+1}} dz_1\wedge dz_2
\]
and then take the residue $\mbox{Res}_{z_1=z_2}f(z_1)dz_1/(z_1-z_2)^{j+1}$.)

Further integrating w.r.t. $z_2$ we obtain the Lie bracket formula
 \begin{eqnarray}
 \label{bracke of lie v - tr inv}
  \mbox{Lie}(V)\otimes \mbox{Lie}V&\rightarrow&  \mbox{Lie}(V)\\
  \overline{f(z_1) a(z_1) dz_1}\otimes \overline{g(z_2) b(z_2) dz_2}&\mapsto&\sum_{j=0}^\infty \frac{1}{j!}\overline{f^{(j)}(z_2)g(z_2)(a_{(j)}b)(z_2) dz_2}\nonumber.
  \end{eqnarray}
  
  To see what this really means denote by $v_{(n)}=\overline{z^n v dz}$ for any $v\in V$. The fact that (\ref{bracke of lie v - tr inv}) defines a Lie algebra structure and  (\ref{acti of lie v - tr inv}) this Lie algebra action means the following identity
  \begin{equation}
   \label{borch-comm-in-n-ope-2}
 [a_{(n)},b_{(m)}]c=\sum_{s=0}^\infty{n\choose s}(a_{(s)}b)_{(n+m-s)}c,\; \forall a,b,c\in V,
 \end{equation}
 which is familiar to us and was proved in sect.~\ref{expl-form-lie-star} only if $m,n\geq 0$ from a slightly different perspective.

  \subsubsection{ Example.} 
  \label{example of 2 inf dim lie alg} 
  In the case of the coisson algebra $J_\infty\BC[T^*\Pi TX/S]\otimes\Omega_\BC$, which was defined in sect. \ref{Example q class},
  the Lie algebra $\mbox{Lie}(J_\infty\BC[T^*\Pi TX/S])$ contains a subalgebra with basis $\overline{z^n\partial_{x_i} dz},\overline{z^{m}x_j dz}$,
  $\overline{z^n\partial_{dx_i} dz},\overline{z^{m}dx_j dz}$, $1\leq i\leq N$, $n\in \BZ$, $C=\overline{dz/z}$,  and brackets
  \begin{equation}
  \label{lie alg of free fields f coeff}
  [\overline{z^n\partial_{x_i} dz},\overline{z^{m-1}x_j dz}]=\delta_{n,-m}\delta_{ij}C= [\overline{z^n\partial_{dx_i} dz},\overline{z^{m-1}dx_j dz}],
  \end{equation}
  all the other feasible pairwise brackets being 0. This follows from (\ref{bracke of lie v - tr inv}) applied to  the formulas recorded in sect.\ref{Example q class}.
  
  Let us temporarily denote this algebra by $\Gamma_N$.
  
  In the case of $J_\infty\BC[T^*X/S]$, we consider the subalgebra of  $\Gamma_N$ generated by all of the  even basis elements, $\overline{z^n\partial_{x_i} dz},\overline{z^{m}x_j dz}$; denote it by 
  $\Gamma^{ev}_N$.

 \subsection{ }
 \label{subs on borch ident2}
 An obvious embedding $M\boxtimes M\hookrightarrow j_*^{(2)}j^{(2)*} M\boxtimes M$ means that a chiral bracket 
 $j_*^{(2)}j^{(2)*} M\boxtimes M\rightarrow \Delta^{(2)}_*M$  determines a Lie* bracket $M\boxtimes M\rightarrow \Delta^{(2)}_*M$, by restriction. Each chiral algebra is therefore a Lie* algebra, and so formulas such as  (\ref{acti of lie v - tr inv} -- \ref{borch-comm-in-n-ope-2}) apply, in particular, to a chiral algebra. We will now write down an identity, which generalizes (\ref{borch-comm-in-n-ope-2}) and follows from a full-fledged chiral Jacobi identity.
  
 We let
  $[.,.]\in P^{ch}_2(\{V\otimes\Omega_C, V\otimes\Omega_C\}, V\otimes\Omega_C)$ be a chiral bracket,  consider the composition $[.,[.,.]]$, evaluate it on
  $a(z_1)\otimes b(z_2)\otimes c(z_3) (z_1-z_3)^n(z_2-z_3)^m(z_1-z_2)^l)$, and then ``Jacobiate'' it. Denote the result by 
  $Jac(a(z_1)\otimes b(z_2)\otimes c(z_3) (z_1-z_3)^n(z_2-z_3)^m(z_1-z_2)^l)$.  The Jacobi identity says 
  \[
  Jac(a(z_1)\otimes b(z_2)\otimes c(z_3) (z_1-z_3)^n(z_2-z_3)^m(z_1-z_2)^l)=0.
  \]
    Writing this down  as it is is an onerous task, similar to but harder than the computation that proved (\ref{borch-comm-in-n-ope}). We will restrict ourselves to this equality written down modulo $\partial_{z_1}, \partial_{z_2}$, that is, 
  \[
  Jac(a(z_1)\otimes b(z_2)\otimes c(z_3) (z_1-z_3)^n(z_2-z_3)^m(z_1-z_2)^l)\equiv 0 \mbox{ mod } \partial_{z_1}, \partial_{z_2}.
  \]
  By definition, this is a statement to the effect that a certain element of $V\otimes\Omega_C$ (in fact of $V$,  in translation-invariant setting) that depends on $a,b,c\in V$ and $n,m,l\in\BZ$ is 0. It is easiest to organize the computation of this element as a repeated integration. Using the $[a,[b,c]]-[b,[a,c]]=[[a,b],c]$ format of the Jacobiator the desired identity becomes, in the notation of (\ref{what ope looks like in tr inva},\ref{chiral bracket formula tr inva})
  
   \begin{eqnarray*}
  \iint_{|z_1-z_3|>|z_2-z_3|}dz_1\wedge dz_2((z_1-z_3)^n(z_2-z_3)^m(z_1-z_2)^{-l})a(z_1)\circ(b(z_2)\circ c(z_3))&-&\\
    \iint_{|z_2-z_3|>|z_1-z_3|}dz_1\wedge dz_2((z_1-z_3)^n(z_2-z_3)^m(z_1-z_2)^{-l})b(z_2)\circ(a(z_1)\circ c(z_3))&=&\\
    \label{borch id gen form1}
      \iint_{|z_2-z_3|>|z_1-z_2|}dz_1\wedge dz_2((z_1-z_3)^n(z_2-z_3)^m(z_1-z_2)^{-l})(a(z_1)\circ(b(z_2))\circ c(z_3).
      \end{eqnarray*}
      
      In this formula, the expression under the integration sign indicates the choice of the variables  and of the Laurent series expansion; e.g. in the 1st integral the variables are $z_1-z_3,z_2-z_3$ and 
  \[
  \frac{1}{z_1-z_2}=\sum_{j=0}^\infty \frac{(z_2-z_3)^j}{(z_1-z_3)^{j+1}};
  \]
  in the 2nd,  the variables are $z_1-z_3,z_2-z_3$ , but
  \[
  \frac{1}{z_1-z_2}=-\sum_{j=0}^\infty \frac{(z_1-z_3)^j}{(z_2-z_3)^{j+1}};
  \]

  and in the last, the variables are $z_1-z_2,z_2-z_3$ and
  \[
  \frac{1}{z_1-z_3}=-\sum_{j=0}^\infty \frac{(z_1-z_2)^j}{(z_2-z_3)^{j+1}}.
  \]
  Integration means computing the residue, that is, the coefficient of $dz_1\wedge dz_2/(z_1-z_3)(z_2-z_3)$ for the first two terms, and of 
  $dz_1\wedge dz_2/(z_1-z_2)(z_2-z_3)$ for the last one.

 It is convenient at this point to do what is usually done -- introduce a field:
 \begin{equation}
 \label{def of vert fie}
 V\ni a\mapsto Y(a,z)\stackrel{\mbox{def}}{=} \sum_{m\in\BZ}a_{(m)}z^{-m-1}\in\mbox{End}_\BC(V)((z,z^{-1}))
 \end{equation}
 and rewrite the last identity as follows:

  \begin{eqnarray*}
  \iint_{|z_1-z_3|>|z_2-z_3|}dz_1\wedge dz_2((z_1-z_3)^n(z_2-z_3)^m(z_1-z_2)^{-l})Y(a,z_1-z_3)Y(b,z_2-z_3)c)&-&\\
    \iint_{|z_2-z_3|>|z_1-z_3|}dz_1\wedge dz_2((z_1-z_3)^n(z_2-z_3)^m(z_1-z_2)^{-l})Y(b,z_2-z_3)Y(a,z_1-z_3)c)&=&\\
    \label{borch id gen form}
      \iint_{|z_2-z_3|>|z_1-z_2|}dz_1\wedge dz_2((z_1-z_3)^n(z_2-z_3)^m(z_1-z_2)^{-l})Y(Y(a,z_1-z_2)b,z_2-z_3)c).
      \end{eqnarray*}
    
  The case $l=0$ of this identity gives the {\em Borcherds commutator} formula:
  
   \begin{equation}
 \label{borch-comm-in-n-ope ver versi}
 [a_{(n)},b_{(m)}]c=\sum_{s=0}^\infty{n\choose s}(a_{(s)}b)_{(n+m-s)}c, \forall m,n\in\BZ,
 \end{equation}
 familiar to us, (\ref{borch-comm-in-n-ope} or \ref{borch-comm-in-n-ope-2}).
 
 The case $n=0$, $l>0$ gives us the {\em normal ordering formula}
 \begin{equation}
 \label{norm orde from borche}
 Y(a_{(-l)}b,z)=\frac{1}{(l-1)!}:Y(a,z)^{(l-1)}Y(b,z):,
 \end{equation}
 where we let
 \begin{eqnarray*}
 Y(a,z)_+&=&\sum_{m=0}^{+\infty}a_{(-m-1)}z^{m},  \;Y(a,z)_-=\sum_{m=-1}^{-\infty}a_{(-m-1)}z^{m}, \\
 :Y(a,z)Y(b,z):&=&Y(a,z)_{+}Y(b,z)+Y(b,z)Y(a,z)_-.
 \end{eqnarray*}
 For example, (\ref{norm orde from borche}) asserts
 \begin{equation}
 \label{norm-orde-formu}
 (a_{(-1)}b)_{(n)}c=\sum_{j=0}^{\infty} a_{(-1-j)}b_{(n+j)}c+ b_{(n-j-1)}a_{(j)}c,
 \end{equation}
which formula is easily seen to be a quantization of (\ref{leib-5}); some details will be supplied in sect.~\ref{Filtration and quasiclassical limit.}.

The case where $l\leq 0$ readily reduces to that of $l=0$, because then $(z_1-z_2)^{-l}$ is a polynomial in $z_1-z_3$ and $z_2-z_3$.

 In fact, one can  prove, see \cite{K}, sect. 4.8, that the identities (\ref{borch-comm-in-n-ope ver versi}, \ref{norm orde from borche}) are essentially equivalent to the definition of a vertex algebra (a.k.a. translation-invariant chiral algebra).
 
 \subsection{ }
 \label{quant over c-n}
 We would like to quantize the sheaf of coisson algebras $\CO_{T^*\Pi TX/S}$, see sect. \ref{Example q class}. It can be done locally, where $X/S$ is an affine variety over $S$ whose tangent sheaf has an abelian basis.  In the case  $X/S=\BC^N_S/S$,  $\BC^N_S=\BC^N\times S$, a quantization, that is, an appropriate vertex algebra has been always known.  
 
 \subsubsection{Construction. }
 \label{constr of quan on aff space}

 Consider the Lie $\CO_S$-algebra 
 $\Gamma_N$
 introduced in sect.~\ref{example of 2 inf dim lie alg}. It  has 2 complementary  abelian Lie subalgebras,  $\Gamma_{N,+}$, $\Gamma_{N,-}$, the former being generated, over $\CO_S$, by
$ \overline{z^n\partial_{x_i} dz}$,  $ \overline{z^n\partial_{dx_i} dz}$, $ \overline{z^{n}x_idz}$, $ \overline{z^{n}dx_i dz}$, $C$ with $n\geq 0$, the latter  by $ \overline{z^n\partial_{x_i} dz}$,  $ \overline{z^n\partial_{dx_i} dz}$, $ \overline{z^{n}x_idz}$, $ \overline{z^{n}dx_i dz}$ with $n<0$. Denote by
\begin{equation}
\label{def of omega ch c n}
\Omega^{ch}(\BC^N_S/S)=\mbox{Ind}_{\Gamma_{N,+}}^{\Gamma_{N}}\BC_1,
\end{equation}
 the  $\Gamma_N$-module induced from $\BC_1$, the 1-dimensional representation of $\Gamma_{N,+}$, where $C$ acts as multiplication by 1, and the rest of the indicated basis elements act trivially.
 
 As a $\Gamma_{N,-}$-module, $\Omega^{ch}(\BC_S^N/S)$ is naturally identified with the universal enveloping algebra $U(\Gamma_{N,-})$. The latter is in fact a polynomial ring, and as such it is naturally identified with the jet algebra $J_\infty\BC[T^*\Pi T\BC_S^N/S]$, see sect. \ref{Example q class}. Namely, we have a composite isomorphism 
 
 \begin{equation}
 \label{iden modul unive enve jetss}
 \Omega^{ch}(\BC_S^N/S)\stackrel{\sim}{\rightarrow}U(\Gamma_{N,-})\stackrel{\sim}{\rightarrow}J_\infty\BC[T^*\Pi T\BC_S^N/S]=\CO_S[ x_i^{(n)},\partial_{x_i}^{(n)},  dx_i^{(n)},\partial_{dx_i}^{(n)};\; n\geq 0],
 \end{equation}
 which acts on generators is as follows:
 \begin{eqnarray*}
\overline{z^{-n-1}x_idz}&\mapsto &x_i^{(n)},\\
\overline{z^{-n-1}\partial_{x_i}dz}&\mapsto &\partial_{x_i}^{(n)},\\
\overline{z^{-n-1}dx_idz}&\mapsto &dx_i^{(n)},\\
\overline{z^{-n-1}\partial_{dx_i}dz}&\mapsto &\partial_{dx_i}^{(n)}
 \end{eqnarray*}
  Elements on the right make sense only if $n\geq 0$; elements on the left,
 $ \overline{z^{-n-1}x_idz}$, $\overline{z^{-n-1}\partial_{x_i} dz}$,
 $ \overline{z^{-n-1}dx_idz}$, $\overline{z^{-n-1}\partial_{dx_i} dz}$, make sense for all $n$, positive and negative. Reflecting this, and in keeping with the standard vertex algebra notation, from now on we will  be always writing
 \begin{eqnarray}
 x_{i,(n)}&\mbox{ for }&  \overline{z^{n}x_idz},\nonumber\\ 
\partial_{ x_{i,(n)}}&\mbox{ for }&  \overline{z^{n}\partial_{x_i}dz},\nonumber\\
 & & \label{some_identifications}\\
dx_{i,(n)}&\mbox{ for }&  \overline{z^{n}dx_idz},\nonumber\\
\partial_{ dx_{i,(n)}}&\mbox{ for }&  \overline{z^{n}\partial_{dx_i}dz}.\nonumber
\end{eqnarray}

 The space $\Omega^{ch}(\BC_S^N/S)\otimes\CO_\BC$ carries a chiral algebra structure defined OPE style as follows. 
 
 The unit, $\bun\in\Omega^{ch}(\BC_S^N/S)$, is $1\in\BC_1$.
 
 $\Omega^{ch}(\BC_S^N/S)$ inherits a derivation $\partial$ from
 $\CO_S[ x_i^{(n)},\partial_{x_i}^{(n)},  dx_i^{(n)},\partial_{dx_i}^{(n)};\; n\geq 0]$. This makes $\Omega^{ch}(\BC_S^N/S)\otimes\CO_\BC$ into a translation-invariant $D_\BC$-module.

 We have $\partial(\bun)=0$ and $\partial(a_{(-n-1)})=(n+1)a_{(-n-2)}$, $n\geq 0$, for $a=x_i, \partial_{x_i}, dx_i, \partial_{dx_i}$; cf. 
 (\ref{action of d on jet}, \ref{expl form lin over d all n})

 To each $a\in\Omega^{ch}(\BC_S^N/S)$ we attach a field 
 \[
 Y(a,z)=\sum_{n\in\BZ}a_{(n)}z^{-n-1}.
 \]
 
 Since we have identified $\Omega^{ch}(\BC_S^N/S)$ with a polynomial ring, the field $Y(a,z)$ will be defined by induction on the degree of $a$. In degree 0,
\[
Y(\bun,z)=\mbox{Id};
\]
in degree 1,
\[
Y(a_{(-k-1)}\bun,z)=Y(a_{(-k-1)},z)= \frac{1}{k!}(\sum_{n\in\BZ}a_{(n)}z^{-n-1})^{(k)} \mbox{ for } a=x_i, \partial_{x_i}, dx_i, \partial_{dx_i};
\]
and more generally,
  \[
 Y( a _{(-k-1)}P,z)=:\frac{1}{k!}Y(a_{(-1)},z)^{(k)}Y(P,z):  \mbox{ for } a=x_i, \partial_{x_i}, dx_i, \partial_{dx_i}, P\in \Omega^{ch}(\BC^N/S).
 \]
which is imposed upon us by (\ref{norm orde from borche}).

This allows us to define an OPE 
 \begin{eqnarray}
\circ:\; (\Omega^{ch}(\BC_S^N/S)\otimes\CO_\BC)\boxtimes (\Omega^{ch}(\BC_S^N/S)\otimes\CO_\BC)&\rightarrow&\widetilde{ \Delta}^{(2)}(\Omega^{ch}(\BC_S^N/S)\otimes\CO_\BC),\nonumber\\
\label{how we defi ope}
a(z_1)\circ b(z_2)&=& Y(a,z_1-z_2)b(z_2),
\end{eqnarray}
for $a,b\in \Omega^{ch}(\BC_S^N/S)\otimes 1$.

\subsubsection{Notational Remark}
\label{Notational Remarks}
The notation $a_{(n)}$,  $a\in \Omega^{ch}(\BC_S^N/S)$,  $n\in\BZ$,  is overburdened and arises in two potentially different contexts. Firstly, as a coefficient of the field
$Y(a,z)=\sum_{n\in\BZ}a_{(n)}z^{-n-1}$, in which case $a_{(n)}$ is a linear transformation of $\Omega^{ch}(\BC_S^N/S)$. Secondly, as a label of a `coordinate function', i.e., of one of the fixed generators of $\Omega^{ch}(\BC_S^N/S)$, such as $\partial_{x_i,(n)}$,  $x_{i,(n)}$, $\partial_{dx_i,(n)}$,
$dx_{i,(n)}$, which happens only if $n<0$, see identifications (\ref{some_identifications}). Therefore, a potential notational clash might only occur in the case of $a_{(-n-1)}$, $n\geq 0$, and $a$ one of those 4 symbols, in which case, luckily, no ambiguity arises: by construction, the operator $a_{(-n-1)}: \Omega^{ch}(\BC_S^N/S)\rightarrow \Omega^{ch}(\BC_S^N/S)$ is in fact multiplication by $a_{(-n-1)}$, so that the omnipresent vertex algebra notation $a_{(-n-1)}b$ means exactly ``$b$ multiplied by $a_{(-n-1)}$," which by the way coincides with the jet space notation $a^{(n)}b$, cf. (\ref{iden modul unive enve jetss}).

\subsubsection{Formulas.}
\label{Associativity. Wick's formula.}
OPE (\ref{how we defi ope}) is algorithmically computable.  For example,
\begin{equation}
\label{appl to one}
a(z_1)\circ\bun(z_2)=\sum_{j=0}^\infty \frac{(z_1-z_2)^j}{j!}\partial^j a(z_2)=\exp{((z_1-z_2)\partial)}a(z_2),
\end{equation}
and 
\begin{equation}
\label{ one appl to }
\bun(z_1)\circ a(z_2)=a(z_2),
\end{equation}
which, by the way, are valid in any vertex algebra.

The Lie bracket in $\Gamma_N$ implies
\begin{equation}
\label{contractions1}
[Y(\partial_{x_i,(-k-1)},z_1-z_2)_-, x_{j,(-l-1)}]=\delta_{ij}{l+k\choose k}\frac{(-1)^k}{(z_1-z_2)^{l+k+1}},
\end{equation}
 \begin{equation}
\label{contractions2}
[Y(x_{i,(-k-1)},z_1-z_2)_-,\partial_{ x_j,(-l-1)}]=-\delta_{ij}{l+k\choose k}\frac{(-1)^k}{(z_1-z_2)^{l+k+1}},
\end{equation}
\begin{equation}
\label{contractions3}
[Y(\partial_{dx_i,(-k-1)},z_1-z_2)_-, dx_{j,(-l-1)}]=\delta_{ij}{l+k\choose k}\frac{(-1)^k}{(z_1-z_2)^{l+k+1}},
\end{equation}
 \begin{equation}
 \label{contractions4}
[Y(dx_{i,(-k-1)},z_1-z_2)_-,\partial_{ dx_j,(-l-1)}]=\delta_{ij}{l+k\choose k}\frac{(-1)^k}{(z_1-z_2)^{l+k+1}}.
\end{equation}
Since $Y(a,z-w)_-\bun=0$, by construction, we can compute the OPE in its entirety by rewriting (\ref{how we defi ope}) as follows:
\begin{equation*}
a(z_1)\circ b(z_2)= (Y(a,z_1-z_2)_+b)(z_2)+(Y(a,z_1-z_2)_-b)(z_2),
\end{equation*}
and then moving $Y(a,z-w)_-$ all the way to the right repeatedly using formulas (\ref{contractions1}--\ref{contractions4}).  For example,
\begin{eqnarray*}
\partial_{x_i,(-1)}(z_1)\circ x_{j,(-l-1)}(z_2)&=&\delta_{ij}\frac{\bun(z_2)}{(z_1-z_2)^{l+1}}+Y(\partial_{x_i,(-1)}, z_1-z_2)_+x_{j,(-l-1)}(z_2)\nonumber\\
\label{ope d vs x}
&=&\delta_{ij}\frac{\bun(z_2)}{(z_1-z_2)^{l+1}}+\exp{((z_1-z_2)\partial})(\partial_{x_i,(-1)}))x_{j,(-l-1)}(z_2).
\end{eqnarray*}
Therefore, if   $x_{j,(-l-1)}$ is replaced with an arbitrary polynomial $b\in\Omega^{ch}(\BC_S^N/S)$, then we obtain

\begin{eqnarray}
\partial_{x_i,(-1)}(z_1)\circ b(z_2)
&=&\sum_{l=0}^{+\infty}\frac{1}{(z_1-z_2)^{l+1}}\frac{\partial b}{\partial x_{i,(-l-1)}}(z_2)+\exp{((z_1-z_2)\partial})(\partial_{x_i,(-1)}))b(z_2)\nonumber\\
\label{ope d vs x 2}
&=&\sum_{l=0}^{+\infty}\frac{1}{(z_1-z_2)^{l+1}}\frac{\partial b}{\partial x_{i,(-l-1)}}(z_2)+\sum_{l=0}^{+\infty}(z_1-z_2)^l\partial_{x_i,(-l-1)}b(z_2).
\end{eqnarray}

The gist of these well-known formulas is the following observation, simple but crucial for what follows: the vertex algebra $\Omega^{ch}(\BC_S^N/S)$ is identified with a commutative associative ring, and the coefficients in front of powers of $(z_1-z_2)$ in $a(z_1)\circ b(z_2)$ are the values of certain bi-differential operators with constant coefficients evaluated on $a\otimes b$. For example,   (\ref{ope d vs x 2}) computes the restriction of this bi-differential operator to the set $\partial_{x_i,(-1)}\otimes \Omega^{ch}(\BC_S^N/S)\approx \Omega^{ch}(\BC_S^N/S)$ to the effect that
 the coefficient  in front of $(z_1-z_2)^{-l-1}$, $l\geq 0$, equals the order 1 differential operator $\partial/\partial x_{i,(-l-1)}$,
  the one in front of $(z_1-z_2)^{l}$, $l\geq 0$, is multiplication by the `function' $\partial_{x_i,(-l-1)}$, that is, an order 0 differential operator; see  sect.~\ref{Notational Remarks} for an elucidation of the symbol $\partial_{x_i,(-l-1)}b(z_2)$.

This process can be iterated, and the result is the following ``explicit'' formula for the OPE. For the purpose of writing down this formula, let 
$I$ stand for the set of linear generators of $\Omega^{ch}(\BC_S^N/S)$ regarded as a polynomial ring -- it is at this point that a choice of a basis
$\{x_i,\partial_{x_i}\}$ becomes essential. Denote by 
\begin{equation*}
\langle a(z_1),b(z_2)\rangle=[Y(a,z_1-z_2)_-,b],
\end{equation*}
the expressions computed in (\ref{contractions1}-\ref{contractions4}); they are called {\em contractions}. We have a somewhat awful formula:

\begin{equation}
\label{expl_form_ope}
a(z_1)\circ b(z_2)= m(\exp{((z_1-z_2)\partial\otimes\mbox{Id})}\exp(\sum_{(u,v)\in I\times I}\langle u(z_1),v(z_2)\rangle\frac{\partial}{\partial u}\otimes \frac{\partial}{\partial v} )(a\otimes b)(z_2)),  
\end{equation}
$m$ on the left standing for the product on $\Omega^{ch}(\BC_S^N/S)$ as an associative ring.
This formula is nothing but a version of the familiar procedure going under the name of Wick's theorem, but it does prove
\subsubsection{ Lemma.}
\label{state ope is defined by diff ope}
{\em There is a collection of bi-differential  operators with constant coefficients on $\Omega^{ch}(\BC_S^N/S)$, $\{P_n,n\in\BZ\}$, such that}
\[
a(z_1)\circ b(z_2)=\sum_{n\in\BZ}(z_1-z_2)^nP_n(a,b)(z_2) .
\]

\subsubsection{ } 
\label{def etale over cn}
The significance of this observation is that thus defined vertex algebra structure on $\Omega^{ch}(\BC_S^N/S)$ can be localized. Namely, suppose
an affine smooth variety $X/S$ is \'etale over $ \BC^N/S$.  We then define 
$\Omega^{ch}_{X/S}=\CO_{X/S}\otimes_{\BC[\BC_S^N/S]}\Omega^{ch}(\BC_S^N/S)$. The differential nature of OPE (\ref{expl_form_ope}) makes it possible to extend it to $\Omega^{ch}_{X/S}$.

{\bf Definition.}{\em  For $X/S$ \'etale over $\BC_S^N/S$, define $\Omega^{ch}_{X/S}$ to be a vertex algebra with OPE (\ref{expl_form_ope}).}

\bigskip

The same construction works in the purely even case and gives us a CDO $D^{ch}_{X/S}$. Its original definition, \cite{MSV}, was slightly different.

\subsection{ Sample computations, filtration, quasiclassical limit.}
\label{ Filtration, quasiclassical limit, sample computations}

\subsubsection{ Sample computations. }
\label{ Sample computations. }
$\Omega^{ch}(X/S)$ is identified, by construction, with the ring $J_\infty\BC[T^*\Pi TX/S]= \BC[X/S][x_{i}^{(n+1)}, \partial_{x_{i}}^{(n)},
dx_{i}^{(n)}, \partial_{dx_{i}}^{(n)}; n\geq 0]$. (Note the shift $n\mapsto n+1$ in $x_i^{(n+1)}$, $x_i^{(0)}=x_i$ being absorbed by $ \BC[X/S]$.)  In other words, each element of the latter ring is a linear combination of monomials
\[
fa_1^{(n_1)}a_2^{(n_2)}\cdots a_k^{(n_k)}, \; f\in \BC[X/S],\; \mbox{ each } a_j= x_{i}, \partial_{x_{i}},
dx_{i}, \partial_{dx_{i}}.
\]
The corresponding element of $\Omega^{ch}(X/S)$ is
\[
a_{1,(-n_1-1)}a_{2, (-n_2-1)}\cdots a_{k,(-n_k-1)}f,
\]
written in this order, see sect.~\ref{Notational Remarks} for an elucidation.  For example, (\ref{expl_form_ope}) will show
\[
\partial_{x_i,(-1)}f =f\partial_{x_i}\mbox{ but }f_{(-1)}(\partial_{x_i,(-1)})=f\partial_{x_i}-\partial(\frac{\partial f}{\partial x_i}).
\]
For this reason, we will be always  writing $\partial_{x_i,(-1)}f$ if we want to indicate the `vector field' $f^i\partial_{x_i}$.

Continuing in the same vein, (\ref{expl_form_ope}) implies
\begin{equation}
\label{ope v field fnctn}
\partial_{x_i,(-1)}f^i(z_1)\circ g(z_2)=\frac{1}{z_1-z_2}f^i\frac{\partial g}{\partial x_i}(z_2)+\cdots
\end{equation}
$\cdots$ indicating regular terms. This result has an obvious geometric interpretation as it reproduces the action of a vector field on a function.

Similarly,
\begin{eqnarray}
\label{ope v field vs v field}
& &\partial_{x_i,(-1)}f^i(z_1)\circ \partial_{x_j,(-1)}g^j(z_2)=\\
\nonumber
& &\frac{1}{(z_1-z_2)^2}\left(-\frac{\partial f^i}{\partial x_j}\frac{\partial g^j}{\partial x_i}\right)(z_2)+\\
\nonumber
& &\frac{1}{z_1-z_2}\left(\partial_{x_j, (-1)}f^i\frac{\partial g^j}{\partial x_i}-\partial_{x_i, (-1)}g^j\frac{\partial f^i}{\partial x_j}-\partial(\frac{\partial f^i}
{\partial x_j})\frac{\partial g^j}{\partial x_i} \right)(z_2)+ \cdots
\end{eqnarray}
\begin{sloppypar}
This result contains two terms that have no obvious geometric meaning, $-(\partial f^i/\partial x_j)(\partial g^j/\partial x_i)$ and 
$-(\partial(\partial f^i/\partial x_j))(\partial g^j/\partial x_i)$, both stemming from the double contractions in (\ref{expl_form_ope}). We will see in sect.~\ref{Calc 2.} how this anomaly disappears in the chiral de Rham situation.
\end{sloppypar}

\subsubsection{Filtration and quasiclassical limit.}
\label{Filtration and quasiclassical limit.}
$\Omega^{ch}(X/S)$ carries an increasing filtration $\{\CF^i\}$, where $\CF^i$ is spanned by monomials containing at most $i$ letters $\partial_\bullet$,
cf. sect. \ref{ Sample computations. }. This filtration is preserved by the OPE and so the corresponding graded object,
$\mbox{Gr}_\CF \Omega^{ch}(X/S)$, is also a vertex algebra. Formula (\ref{expl_form_ope}) shows that the corresponding OPE has no poles, which  
makes $\mbox{Gr}_{\CF^\bullet} \Omega^{ch}(X/S)$ equipped with product $_{(-1)}$ into a commutative associative algebra with derivation $\partial$, obviously isomorphic to $J_\infty\BC[T^*\Pi TX/S]$.   Furthermore, picking  in the standard manner the subleading term of the OPE one obtains a coisson bracket.
A moment's thought will show that doing so amounts to keeping in (\ref{expl_form_ope}) only those terms where the operators of the type $\partial/\partial a\otimes\partial/\partial b$ are applied  once; this clearly makes (\ref{expl_form_ope}) into (\ref{bracket in general even}). We conclude that
$\Omega^{ch}(X/S)$ is indeed a quantization of $J_\infty\BC[T^*\Pi TX/S]$.

As an example, a glance at (\ref{ope v field vs v field}) will show that the passage to the graded object removes the anomalous term, $\partial(\partial f^i/
\partial x_j)\partial g^/\partial x_i$, and makes sure that
\[
(\partial_{x_i,(-1)}f^i)_{(0)}( \partial_{x_j,(-1)}g^j(z_2))=\{f^i\partial_{x_i},g^j\partial_{x_j}\},
\]
the latter being the ordinary Lie bracket of vector fields (in the present context arising as the Poisson bracket of fiberwise-linear functions on $T^*X/S$.)

As another example,  note that picking the subleading term in formula (\ref{norm-orde-formu}) with $n\geq 0$, which amounts to keeping only those $b_{(n-j-1)}a_{(j)}$ where $j\geq n$, reproduces its quasiclassical limit, formula (\ref{leib-5}).

\section{chiral calculus}
\label{chiral calculus}

\subsection{Calc 2.}
\label{Calc 2.}
It is a striking feature of the chiral de Rham complex \cite{MSV} that many of the classical formulas reproduced in sect. \ref{Calc1.} carry over to
$\Omega^{ch}(X/S)$.  Let us explain this.

 Similarly to (\ref{ordinary picard lie}),  the term $\CF^1\Omega^{ch}(X/S)$  fits the exact sequence
\begin{equation}
\label{chiral picard lie}
0\rightarrow J_\infty\Omega^\bullet(X/S)\rightarrow \CF^1\Omega^{ch}(X/S)\rightarrow J_\infty\CT_{\Pi TX/S}\rightarrow 0
\end{equation}
and can be called a Picard-Lie {\em chiral} algebroid, \cite{BD}; the Lie* algebra $J_\infty\CT_{\Pi TX/S}$ was introduced in sect. \ref{jets in lie alg current alg}. As a practical matter, and by the definition of the filtration $\{\CF^i\}$, sect.~\ref{Filtration and quasiclassical limit.},   $\CF^1\Omega^{ch}(X/S)$ is the span of terms that contain at most one symbol $\partial_\bullet$, and $J_\infty\Omega^\bullet(X/S)=\CF^0\Omega^{ch}(X/S)$ is the span of those that contain no such symbol; but conceptually it carries the following structure:

$\bullet$  $\Omega^{ch}(X/S)$ is a vertex algebra and $J_\infty\Omega^\bullet(X/S)$ is its commutative subalgebra;

$\bullet$  being a vertex algebra, $\Omega^{ch}(X/S)$ is a Lie* algebra and  $\CF^1\Omega^{ch}(X/S)$ is its Lie* subalgebra of which $J_\infty\Omega^\bullet(X/S)$ is an abelian ideal and  $J_\infty\CT_{\Pi TX/S}$ is a Lie* algebra naturally isomorphic to the quotient $\CF^1\Omega^{ch}(X/S)/J_\infty\Omega^\bullet(X/S)$;

$\bullet$  $J_\infty\CT_{\Pi TX/S}$ is properly called a Lie* algebroid because it is a Lie* algebra, which acts on $J_\infty\Omega^\bullet(X/S)$  by derivations,  and a module over $J_\infty\Omega^\bullet(X/S)$ and these 2 structures are compatible in that the action map $J_\infty\Omega^\bullet(X/S)\otimes^! J_\infty\CT_{\Pi TX/S}\rightarrow J_\infty\CT_{\Pi TX/S}$ is a morphism of $J_\infty\CT_{\Pi TX/S}$-modules; this is similar to Definition \ref{def of coiss alg}.

Unlike (\ref{ordinary picard lie}), the sequence (\ref{chiral picard lie}) does not split. Nevertheless, (\ref{formula for lie br on diffe forms}, \ref{formula for contr with a vect fie on diffe forms}) can be ``lifted'' to $\Omega^{ch}_{X/S}$ as follows.

For $f_i\partial_{x_i}\in \CT_{X/S}$, define
\begin{equation}
\label{chiral_vers_of_Lie_deriv}
\mbox{Lie}^{ch}_{f_i\partial_{x_i}}=\partial_{x_i,(-1)}f_i-\partial_{dx_i,(-1)}\frac{\partial f_i}{\partial x_\alpha}dx_\alpha\in \CF^1\Omega^{ch}(X/S);
\end{equation}
this can be thought of as a `chiralization', and a naive one at that, of (\ref{expl formula for lie derib and contra ordin}).

Computing as in sect. \ref{ Sample computations. }, we obtain the OPE

\begin{eqnarray}
\label{ope v field vs v field in deram}
\nonumber
& &\mbox{Lie}^{ch}_{f^i\partial_{x_i}}(z_1)\circ \mbox{Lie}^{ch}_{g^j\partial_{x_j}}(z_2)=\\
\nonumber
& &\frac{1}{(z_1-z_2)^2}\left(-\frac{\partial f^i}{\partial x_j}\frac{\partial g^j}{\partial x_i}+\frac{\partial f^\alpha}{\partial x_\beta}\frac{\partial g^\beta}{\partial x_\alpha}\right)(z_2)+
\frac{1}{z_1-z_2}\mbox{Lie}^{ch}_{[f^i\partial_{x_i}g^j\partial_{x_j}]}(z_2)+ \cdots
\end{eqnarray}

This computation should be compared with (\ref{ope v field vs v field}); the reason the coefficient of $1/(z_1-z_2)^2$ vanishes is that now, when using 
(\ref{expl_form_ope}),  there are 2 contributions, one coming from the double contraction of even variables $\partial_{x_i}$ and $x_i$, 
\[
\langle \partial_{x_i,(-1)}(z_1), x_{i,(-1)}(z_2)\rangle\cdot\langle x_{j,(-1)}(z_1), \partial_{x_j,(-1)}(z_2)\rangle
\]
 and another one,  coming from the double contraction of odd variables $\partial_{dx_i}$ and $dx_i$, 
 \[
\langle \partial_{dx_\alpha,(-1)}(z_1), dx_{\alpha,(-1)}(z_2)\rangle\cdot\langle dx_{\beta,(-1)}(z_1), \partial_{dx_\beta,(-1)}(z_2)\rangle
\]
 and they cancel against each other.

Therefore,
\begin{equation}
\label{ope v field vs v field in deram2}
\mbox{Lie}^{ch}_{\xi}(z_1)\circ \mbox{Lie}^{ch}_{\eta}(z_2)=\frac{1}{z_1-z_2}\mbox{Lie}^{ch}_{[\xi,\eta]}(z_2)+\cdots.
\end{equation}

For $f_i\partial_{x_i}\in \CT_{X/S}$, similarly define

\begin{equation}
\label{chiral_versi_of_contract}
\iota^{ch}_{f_i\partial_{x_i}}=\partial_{dx_i,(-1)}f_i,
\end{equation}
a chiralization of the 2nd half of  (\ref{expl formula for lie derib and contra ordin}).

An analogous but considerably easier computation will give

\begin{equation}
\label{ope v field vs contraction with v field in deram2}
\mbox{Lie}^{ch}_{\xi}(z_1)\circ \iota^{ch}_{\eta}(z_2)=\frac{1}{z_1-z_2}\iota^{ch}_{[\xi,\eta]}(z_2)+\cdots.
\end{equation}

We also have a ``chiral version'' of the de Rham differential: 
\[
d^{ch}_{X/S}=\partial_{x_i,(-1)}dx_i,
\]
which satisfies
\begin{equation*}
d_{X/S}(z_1)\circ d_{X/S}(z_2)=0+\cdots,
\end{equation*}
\begin{equation}
\label{d iota equals Lie}
d^{ch}_{X/S}(z_1)\circ \iota^{ch}_\xi(z_2)=\frac{1}{z_1-z_2}\mbox{Lie}^{ch}_\xi(z_2)+\cdots,
\end{equation}
(no double contractions involved), and a little more cumbersome
\begin{equation}
\label{deram diff not quite invar}
d^{ch}_{X/S}(z_1)\circ \mbox{Lie}^{ch}_{f^i\partial_{x_i}}(z_2)=\frac{1}{(z_1-z_2)^2}\left(-\frac{\partial^2f_i}{\partial x_i\partial x_\alpha}dx_\alpha\right)(z_2)+\cdots,
\end{equation}
ellipsis standing for terms regular on the diagonal (no poles involved.) This  implies
\begin{equation}
\label{deram diff 0 is quite}
d^{ch}_{X/S,(0)}\mbox{Lie}^{ch}_{f^i\partial_{x_i}}(z_2)=0.
\end{equation}

It is well known, and follows from (\ref{borch-comm-in-n-ope ver versi}), that for either of the structures we have considered, be it a Lie*, coisson or chiral algebra, if $v\in V$, then $v_{(0)}$ is a derivation of the corresponding operation; therefore (\ref{chiral picard lie}) is a chiral algebroid with differential  $d^{ch}_{X/S,(0)}$.

As $J_\infty\CT_{X/S}$  is a free $\BC[\partial]$-module,  sect. \ref{jets in lie alg current alg}, the maps we have just defined, $\iota^{ch}$ and 
$\mbox{Lie}^{ch}$ uniquely extend to the $\BC[\partial]$-module maps, to be denoted in the same way,
\[
\mbox{Lie}^{ch}: J_\infty\CT_{X/S}\rightarrow \CF^1\Omega^{ch}(X/S);\; \iota^{ch}: J_\infty \Pi\CT_{X/S}\rightarrow \CF^1\Omega^{ch}(X/S).
\]

Formulas (\ref{ope v field vs v field in deram2}, \ref{ope v field vs contraction with v field in deram2}, \ref{d iota equals Lie}) 
 imply the following analogue of Lemma \ref{summ of commu rel on the ordn deram}.

\begin{lem}
\label{summ of commu rel on the ordn deram chiral}
The map
\begin{equation*}
\iota^{ch}\oplus \mbox{Lie}^{ch}:\;J_\infty\Pi\CT_{X/S}\rtimes J_\infty\CT_{X/S}\rightarrow  \CF^1\Omega^{ch}(X/S),
\end{equation*}
defines a morphism of differential  Lie* algebras,
the differential on both sides being $d_{X/S,(0)}$.  (The components $\iota^{ch}$ and  $\mbox{Lie}^{ch}$ are defined in (\ref{chiral_versi_of_contract}) and  (\ref{chiral_vers_of_Lie_deriv}) resp.)
\end{lem}

As we have seen, sect. \ref{subs on borch ident1}, attached to any Lie* algebra  $V$ on a punctured disc is a Lie algebra $\mbox{Lie}(V)$. If $V$ is a differential Lie* algebra, then $\mbox{Lie}(V)$ is also a differential Lie algebra. Formula (\ref{bracke of lie v - tr inv}) shows that $\mbox{Lie}(J_\infty \fg)$ is the loop algebra 
$\fg((t))$, $\mbox{Lie}(J_\infty \fg)\ni \overline{z^n dz\otimes g}\mapsto g\otimes t^n\in\fg((t))$. In our case, 
Lemma \ref{summ of commu rel on the ordn deram chiral} gives us
\begin{cor}
\label{loop alge interpr}
The map $\iota^{ch}\oplus \mbox{Lie}^{ch}$ defines on $ \CF^1\Omega^{ch}(X/S)$ a structure of a differential module over the differential Lie algebra
$\Pi\CT_{X/S}((t))\rtimes \CT_{X/S}((t))$, the differential acting  on the Lie algebra  as follows: $(a\otimes t^n,b\otimes t^m)\mapsto (0,a\otimes t^n)$.
\end{cor}

There is an obvious and natural map $\Omega^\bullet_{X/S}=\BC[\Pi TX/S]\hookrightarrow \Omega^{ch}(X/S)=J_\infty\BC[T^*\Pi TX/S]$, which sends a differential form to itself but now regarded as a function on $\Pi TX/S$. Here is how  classical constructions of sect. \ref{Calc1.} are formally related to their chiralization.

\begin{cor} 
\label{a morp fromn classic to chira}
The map $\Omega^\bullet_{X/S}\hookrightarrow \Omega^{ch}(X/S)$ is a morphism of differential $\Pi\CT_{X/S}\rtimes\CT_{X/S}$-modules with the action on the left defined via Lemma \ref{summ of commu rel on the ordn deram}, and on the right via Corollary \ref{loop alge interpr}: $\iota_\xi$ acts as 
$\iota^{ch}_{\xi, (0)}$ and $\mbox{Lie}_\xi$  as $\mbox{Lie}^{ch}_{\xi, (0)}$.
\end{cor}

\subsection{Calc 3.}
\label{Calc 3.}
We will now chiralize the results of sect. \ref{boryas constr ordin}.

 By construction, the vertex algebra $\Omega^{ch}(X/S)$ has $\CO_S$ as the ground ring. Any commutative associative unital algebra can be regarded as a vertex algebra with 0 derivation, and so we have a vertex algebra morphism $\CO_S\rightarrow \Omega^{ch}(X/S)$. We will now further extend the scalars and introduce $\Omega^\bullet_S\otimes_{\CO_S}\Omega^{ch}(X/S)$, which is also a vertex algebra with the de Rham complex $\Omega^\bullet_S$ as the ground ring. This algebra will be endowed with a differential that can be thought of as a {\em chiral Gauss-Manin connection.}
 
 We continue working strictly locally and assume that $S$ is affine and carries a coordinate system $\{u_i,\partial_{u_i}\}$, $\partial_{u_i}(u_j)=\delta_{ij}$.
 Fix an arbitrary connection on  $X\rightarrow S$, written in coordinates as follows: $\nabla= du_i(\partial_{u_i}+\xi_i)$, $\xi_i\in\CT_{X/S}$.  If we want to emphasize that it defines an action of vector fields on the de Rham complex, we will write $\mbox{Lie}_\nabla=  du_i(\partial_{u_i}+\mbox{Lie}_{\xi_i})$.
 
 We have  the curvature  of $\nabla$:
 \begin{equation}
 \label{def of curv}
 R=\frac{1}{2}[\nabla,\nabla]= \frac{1}{2}du_idu_j(\partial_{u_i}\xi_j-\partial_{u_j}\xi_i+[\xi_i,\xi_j]).
 \end{equation}
 
 Thanks to Corollary \ref{a morp fromn classic to chira}, we know how to lift  these to an operator acting on  $\Omega^\bullet_S\otimes_{\CO_S}\Omega^{ch}(X/S)$:  we define
 \[
 \mbox{Lie}^{ch}_\nabla=  du_i(\partial_{u_i}+\mbox{Lie}^{ch}_{\xi_i})\mbox{ and }
 \]
 \[
 \mbox{Lie}^{ch}_{\nabla,(0)}=  du_i(\partial_{u_i}+\mbox{Lie}^{ch}_{\xi_i,(0)});
 \]
 \[
 \mbox{Lie}^{ch}_R= \frac{1}{2}du_idu_j(\mbox{Lie}^{ch}_{\partial_{u_i}\xi_j}-\mbox{Lie}^{ch}_{\partial_{u_j}\xi_i}+\mbox{Lie}^{ch}_{[\xi_i,\xi_j]})\mbox{ and }
 \]
 \[
 \mbox{Lie}^{ch}_{R,(0)}= \frac{1}{2}du_idu_j(\mbox{Lie}^{ch}_{\partial_{u_i}\xi_{j, (0)}}-\mbox{Lie}^{ch}_{\partial_{u_j}\xi_{i,(0)}}+\mbox{Lie}^{ch}_{[\xi_i,\xi_j],(0)}),
 \]
 and note that thanks to Lemma~\ref{a morp fromn classic to chira}, definition (\ref{def of curv}) implies
 \begin{equation}
 \label{def of curv chir}
  \mbox{Lie}^{ch}_{R,(0)}=\frac{1}{2}[ \mbox{Lie}^{ch}_{\nabla,(0)}, \mbox{Lie}^{ch}_{\nabla,(0)}].
  \end{equation}
  
  Consider $D^{ch}_\nabla\in\mbox{End}_\BC(\Omega^\bullet_S\otimes_{\CO_S}\Omega^{ch}(X/S))$ defined by 
  \begin{equation*}
  D^{ch}_\nabla=\mbox{Lie}^{ch}_{\nabla,(0)}+ d_{X/S,(0)}-\iota^{ch}_{R,(0)}.
  \end{equation*}
  \begin{lem} (cf. Lemma \ref{boryas diff ordina}.) The pair $(\Omega^\bullet_S\otimes_{\CO_S}\Omega^{ch}(X/S),D_\nabla)$ is a complex, that is,
  \label{chiral square is zero}
  \[
  (D^{ch}_\nabla)^2=0.
  \]
  \end{lem}
  
  {\em Proof:}
  \begin{eqnarray*}
  & &2 (D^{ch}_\nabla)^2=[ D^{ch}_\nabla,  D^{ch}_\nabla]=\\
 & & [\mbox{Lie}^{ch}_{\nabla,(0)},\mbox{Lie}^{ch}_{\nabla,(0)}]+[ d_{X/S,(0)}, d_{X/S,(0)}]+[\iota^{ch}_{R,(0)},\iota^{ch}_{R,(0)}]\\
 &+&2[\mbox{Lie}^{ch}_{\nabla,(0)},d_{X/S,(0)}]-2[\mbox{Lie}^{ch}_{\nabla,(0)}, \iota^{ch}_{R,(0)}]-2[d_{X/S,(0)},\iota^{ch}_{R,(0)}]=\\
 & & 2\mbox{Lie}^{ch}_{R,(0)}-2\iota^{ch}_{[\nabla,R],(0)}-2\mbox{Lie}^{ch}_{R,(0)}=0.
 \end{eqnarray*}
  
  A quick guide to this computation: in the last line, the 1st occurrence of $\mbox{Lie}^{ch}_{R,(0)}$ is via $\mbox{Lie}^{ch}_{R,(0)}=
  [\mbox{Lie}^{ch}_{\nabla,(0)},\mbox{Lie}^{ch}_{\nabla,(0)}]$, which is  (\ref{def of curv chir}), the second is via
  $\mbox{Lie}^{ch}_{R,(0)}=[d_{X/S,(0)},\iota^{ch}_{R,(0)}]$, which is a consequence of (\ref{d iota equals Lie}); finally, $\iota^{ch}_{[\nabla,R],(0)}$ appears because it equals $[\mbox{Lie}^{ch}_{\nabla,(0)}, \iota^{ch}_{R,(0)}]$ according to Lemma~\ref{a morp fromn classic to chira} and vanishes because $[\nabla,R]=0$, which is the Bianchi identity. The vanishing of the rest of the summands is more or less clear. $\qed$
  
  We will call the complex $(\Omega^\bullet_S\otimes_{\CO_S}\Omega^{ch}(X/S),D_\nabla)$ the {\it local  chiral 
  connection complex} and denote it $CC^{ch}(X/S, \nabla)$. 
  
   The choice of $\nabla$ is arbitrary, but different such complexes are canonically isomorphic.
  
  \begin{lem} (cf. Lemma \ref{iso of diff connectns ordinary}.)
  \label{iso of diff connectns}
  The map
  \[
  \exp{(\iota^{ch}_{\nabla_2-\nabla_1,(0)})}:(\Omega^\bullet_S\otimes_{\CO_S}\Omega^{ch}(X/S),D_{\nabla_1})\rightarrow
  (\Omega^\bullet_S\otimes_{\CO_S}\Omega^{ch}(X/S),D_{\nabla_2})
  \]
  is an isomorphism of complexes.
  \end{lem}
  
  {\em Sketch of Proof:}
  
  It is clear that the operator $\iota^{ch}_{\nabla_2-\nabla_1,(0)}\in\mbox{End}_\BC(\Omega^\bullet_S\otimes_{\CO_S}\Omega^{ch}(X/S))$ is nilpotent
  and, as we have once mentioned, a derivation of the chiral bracket; hence it defines an automorphism of the vertex algebra.
  What we need to check then is that 
  \[
   \exp{(ad_{\iota^{ch}_{\nabla_2-\nabla_1,(0)}})}(D^{ch}_{\nabla_1})=D^{ch}_{\nabla_2},
   \]
   where $ad_{\iota^{ch}_{\nabla_2-\nabla_1,(0)}}$ means the operation of taking  the bracket, $[\iota^{ch}_{\nabla_2-\nabla_1,(0)},.]$. A direct but little appetizing computation, using the various identities as in the proof of  Lemma \ref{iso of diff connectns}, will show that
   \begin{eqnarray*}
   \exp{(ad_{\iota^{ch}_{\nabla_2-\nabla_1,(0)}})}(D^{ch}_{\nabla_1})&=&D^{ch}_{\nabla_1}+[\iota^{ch}_{\nabla_2-\nabla_1,(0)}, D^{ch}_{\nabla_1}]\\
   &+&
   \frac{1}{2}[\iota^{ch}_{\nabla_2-\nabla_1,(0)},[\iota^{ch}_{\nabla_2-\nabla_1,(0)}, D^{ch}_{\nabla_1}]]=D^{ch}_{\nabla_2}.
   \end{eqnarray*}

   \subsection{Globalization. }
  
   \label{globalization}
    \begin{sloppypar}The above made computations have a number of consequences.
     \end{sloppypar}
   
   (1) Corollary \ref{a morp fromn classic to chira} says that the natural action of the Lie algebra of vector fields on the ordinary de Rham complex lifts to that on the chiral de Rham complex $\Omega^{ch}(X/S)$ if $X/S$ is sufficiently ``small'' (that is, allows a global coordinate system.) A standard argument
   that goes under the Gelfand-Kazhdan formal geometry moniker implies that for each smooth $X/S$ over $S$ there is a  (sheaf of) differential graded vertex algebra $\Omega^{ch}_{X/S}$, along with an embedding of a commutative differential algebra 
   $\Omega^\bullet_{X/S}\hookrightarrow \Omega^{ch}_{X/S}$. The assignment $X/S\mapsto \Omega^{ch}_{X/S}$ is functorial in $X/S$. All of this was originally proved \cite{MSV} by explicitly  writing  down the universal transition functions that define $\Omega^{ch}_{X/S}$ ( only in the absolute setting which, however, hardly makes any difference.) The crucial formula (\ref{ope v field vs v field in deram2}) appeared in \cite{MS}.
   
   (2) Furthermore, we have, keeping to our assumptions on $X/S$, 
   \begin{thm}
   \label{main}
   There is a sheaf of differential vertex algebras on $X/S$ that is locally isomorphic to the one defined in Lemma~\ref{chiral square is zero} with transition functions as in Lemma~\ref{iso of diff connectns}. Different such complexes are canonically isomorphic.
   \end{thm}
   
   {\em Proof} is a repetition, practically verbatim, of the argument from sect.~\ref{Globalization-Ordin.}.  We pick an open cover $\{U_i\}$ as in sect.~\ref{Globalization-Ordin.} and choose a connection $\nabla_i$ over each $p(U_i)\subset S$. By Lemma~\ref{chiral square is zero}, this gives us a complex
   $((\Omega^\bullet_S\otimes_{\CO_S}\Omega^{ch}_{X/S})|_{U_i},D_{\nabla_i})$.  Let us for the sake of brevity re-denote
   \[
   \Omega^{ch}_i=(\Omega^\bullet_S\otimes_{\CO_S}\Omega^{ch}_{X/S})|_{U_i}.
   \]

   On the intersections, we have isomorphisms, Lemma~\ref{iso of diff connectns},
   \begin{equation*}
   G_{ji}\stackrel{\mbox{def}}{=} \exp{(\iota^{ch}_{\nabla_j-\nabla_i,(0)})}:\Omega^{ch}_i|_{U_i\cap U_j}\rightarrow \Omega^{ch}_j|_{U_i\cap U_j}.
 \end{equation*}
 It is clear that
 \[
 G_{ki}=G_{kj}\circ G_{ji}.
 \]
  Therefore, having fixed a cover $\{U_i\}$ and connections $\{\nabla_i\}$, we obtain a new complex, $\Omega^{ch}_{\nabla}$, by tearing the original 
  $\Omega^\bullet_S\otimes_{\CO_S}\Omega^{ch}_{X/S}$ apart and re-gluing the pieces $\{(\Omega_i,\nabla_i)\}$ using  $\{G_{ij}\}$ as transition functions.
  
  Two distinct such complexes, $\Omega^{ch}_{\nabla}$ and $\Omega^{ch}_{\nabla'}$, are canonically isomorphic. To see this, assume, as we can, that both are defined using the same open cover $\{U_i\}$. Then the collection of isomorphisms 
  \begin{equation*}
 \exp{(\iota^{ch}_{\nabla'_i-\nabla_i,(0)})}:\Omega^{ch}_{\nabla}|_{U_i}\rightarrow\Omega^{ch}_{\nabla'}|_{U_i}
  \end{equation*}
  well defines an isomorphism of complexes $\Omega^{ch}_{\nabla}\rightarrow\Omega^{ch}_{\nabla'}$. $\qed$
  
  We will call the sheaf of differential vertex algebras from the above theorem the {\it global chiral connection complex} and denote it $CC^{ch}(X/S)$.

\bigskip\bigskip


\bigskip\bigskip

F.M.: Dept. of Mathematics, University of Southern California, Los Angeles CA 90089, USA

Email: {\tt fmalikov@usc.edu}

\smallskip

V.S.: Institut de Math\'ematiques de Toulouse, Universit\'e Paul Sabatier, 118 Route de Narbonne, 
31062 Toulouse, France. 

Email: 
 {\tt schechtman@math.ups-tlse.fr }
 
 \
 
B.T.: Dept. of Mathematics, Northwestern University, Evanston, Illinois 60208-2730

Email: 
{\tt b-tsygan@northwestern.edu}

\newpage



\end{document}